\documentclass[12pt,fleqn]{scrartcl} 

\usepackage[T1]{fontenc} 
\usepackage[latin1]{inputenc} 

\usepackage{framed}
\usepackage{hyperref}
\usepackage{doi}

\usepackage{amssymb}
\usepackage{amsmath}
\usepackage{bm}

\usepackage{color}
\newtheorem{theorem}{Theorem}%  meant for continuous numbers
\newtheorem{proposition}[theorem]{Proposition}% 
\newtheorem{remark}{Remark}%
\newtheorem{lemma}{Lemma}

\newtheorem{definition}{Definition}%

\newcommand{\BR}{{\mathbb R}}

\newcommand{\BN}{{\mathbb N}}

\newcommand{\D}{\bm D}
\newcommand{\Laplace}{\bm \Delta}
\newcommand{\X}{\bm x}
\newcommand{\XI}{\bm \xi}

\newcommand{\Y}{\bm y}
\newcommand{\E}{{\bm E}}

\newcommand{\cl}{C \kern -0.1em \ell}

\newcommand{\e}{{\bm e}}

\newcommand{\dx}{{\bm d}{\bm x}}
\newcommand{\dy}{{\bm d}{\bm y}}
\newcommand{\dxi}{{\bm d}{\bm \xi}}

\newcommand{\f}{\bm f}
\newcommand{\g}{\bm g}

\renewcommand{\u}{\bm u}
\newcommand{\PSI}{\bm \psi}

\usepackage{soul, cancel, xcolor}

\definecolor{dred}{RGB}{191 0 64}

\usepackage{ulem}
\begin{document}
\title{Paley-Wiener Type Theorems associated to Dirac Operators of Riesz-Feller type}
\date{Swanhild Bernstein\footnote{\textsl{swanhild.bernstein@math.tu-freiberg.de}, Institute of Applied Analysis, TU Bergakademie Freiberg, Pr\"uferstr. 9, D-09599 Freiberg, Germany}, Nelson Faustino\footnote{\textsl{nfaust@ua.pt}, Department of Mathematics, University of Aveiro, Campus Universit\'ario de Santiago, Aveiro, 3810-193, Aveiro, Portugal}}

\maketitle

\begin{abstract}
The classical Paley-Wiener theory consists of two parts. First, the equivalence of the compact support of the Fourier transform of the boundary values to the exponential growth of the extended function, and second, the equivalence to the Bernstein spaces. 
This paper explores Paley-Wiener-type theorems within the framework of hypercomplex variables.
The investigation focuses on a space-fractional version of the Dirac operator $\D_\theta^\alpha$ of order $\alpha$ and skewness $\theta$. The pseudo-differential reformulation of $\D_\theta^\alpha$ in terms of the Riesz derivative $(-\Delta)^{\frac{\alpha}{2}}$ and the so-called \textit{Riesz-Hilbert transform} $\mathcal{H}$, allows for the description of generalized Hardy spaces,
% on the upper and lower half-spaces of $\BR^{n+1}$, $\BR^{n+1}_+$ resp. $\BR^{n+1}_-$, 
using L\'evy-Feller type semigroups generated by $-(-\Delta)^{\frac{\alpha}{2}}$, and the boundary values.
% $\f_\pm=\frac{1}{2}\left(\f\pm \mathcal{H}\f\right)$.

Subsequently, we employ a proof strategy rooted in {\textit real Paley-Wiener methods} to demonstrate that the growth behavior of the sequences of functions $\left(\left(~\D_\theta^\alpha~\right)^k\f_\pm\right)_{k\in \BN_0}$ effectively captures the relationship between the support of the Fourier transform $\widehat{\f}=\mathcal{F}\f$ of the $L^p-$function $\f$, in the case where  $\operatorname{supp}\widehat{\f}\subseteq \overline{B(0,R)}$, and the solutions of Cauchy problems equipped with the space-time operator $\partial_{x_0} + \D_\theta^\alpha$, which are of exponential type $R^\alpha$.

Within the developed framework, introducing a hypercomplex analogue for the Bernstein spaces $B_R^p$ arises naturally, allowing for the meaningful extension of the results by Kou and Qian as well as Franklin, Hogan, and Larkin. 

\end{abstract}


\section{Introduction}\label{sec1}
\subsection{State of art}

Paley-Wiener type theorems, introduced shortly after the celebrated monograph \cite{Paley-Wiener1934} by Paley and Wiener, is a cornerstone approach at the intersection of complex analysis, harmonic analysis, and sampling theory that provides a faithful description of the space of functions or distributions through growth at infinity and the support of their Fourier-like transform.

The classical Paley-Wiener theorem is mainly referred to $L^2-$functions with exponential growth rate as $\vert x+iy\vert\rightarrow \infty$. Let us denote by $\mathbb{C}^+$ the upper complex half-space, which can be identified with $\mathbb{C}^+\cong \mathbb{R}\times (0,+\infty)$. Then, the  Hardy space $H^2(\mathbb{C}^+)$ consists on analytic functions $F$ on $\mathbb{C}^+$ with finite norm
$$ \|F\|_{H^2(\mathbb{C}^+)} := \left( \sup_{y>0}  \int_{\BR} \vert F(x+iy)\vert^2 dx \right)^{\frac{1}{2}} .$$

Specifically, the result obtained in \cite[Theorem X]{Paley-Wiener1934}  (see also~\cite[Chapter III]{SteinWeiss72}) states that the action of the Fourier transform allows us to establish an isometric isomorphism between $H^2(\mathbb{C}^+)$ and $L^2(0,+\infty)$. 

Generalizations of Paley-Wiener type theorems can also be obtained for Bernstein spaces $B_R^p$, with $R>0$ and $1\leq p\leq \infty$. It is known from \cite[p.~98]{Boas54} that $f\in B_{R}^p$ if, and only if, $f$ is an entire function satisfying 
\begin{eqnarray*}
	\left\|~f(\cdot+iy)~\right\|_p \leq \left\|~f ~\right\|_p e^{R \vert y \vert}, &y\in \BR,
\end{eqnarray*} where $ \vert\vert \cdot \vert\vert_p$ denotes the standard $L^p-$norm on $\BR$.

The Bernstein spaces $B_{R}^p$ can be characterized in different ways:
\begin{itemize}
	\item[\textbf (i)] A function $f\in L^p(\mathbb{R})$ belongs to $B_{R}^p$ if, and only if, its distributional Fourier transform has support $[-R , R]$ in the sense of distributions.
	\item[\textbf (ii)] Let $f\! \in\! C^{\infty}(\mathbb{R})$ be such that $f^{(k)}\!\in \!L^p(\mathbb{R})$, for all $k\!\in\! \mathbb{N}_0$ and some $1\leq p\leq \infty$. Then, $f\in B_{R}^p$ if and only if $f$ satisfies the so-called Bernstein's inequality (cf.~\cite{Nik75}, p.~116)
	$$ \vert\vert f^{(k)} \vert\vert_p \leq R ^k \vert\vert f \vert \vert_p, \quad k\in \mathbb{N}_0. $$
\end{itemize}

The subsequent characterization, which establishes the equivalence between Bernstein's inequalities and Bernstein-type spaces, is currently referred to as the {\textit real Paley-Wiener type theorems} (see, e.g., \cite{Andersen04, AndersenJeu10, Andersen14}). This field has been the focus of extensive research by numerous scholars since the seminal contributions of Pesenson \cite{Pesenson88}, Krein-Pesenson \cite{KreinPesenson90}, Bang \cite{Bang90}, and Tuan \cite{Tuan99}. One of the primary reasons for this interest is that it can be universally extended to other Fourier-like transforms, whereas the complex Paley-Wiener approach cannot do so (cf. \cite{TuanZayed02,Andersen06,Pesenson08,PesensonZayed09}).

In the hypercomplex setting, Kou and Qian made the first successful attempt, documented in their seminal paper \cite{KouQian2002}. In particular, they showed that the lack of a generalization for the Phragm\'en-Lindel\"of theorem (cf. \cite{Boas54}) to hypercomplex variables did not prevent them from proving their result. Years later, Franklin, Hogan, and Larkin reformulated Kou and Qian's proofs in terms of Bernstein-type spaces in \cite{FranklinHoganLarkin17}. Recently, Li, Leng, and Fei provided a different approach to real Paley-Wiener spaces in \cite{LiLengFei2019}.

In the papers \cite{FranklinHoganLarkin17, LiLengFei2019}, the authors use a Clifford-type Fourier transform. This approach first appeared in the context of hypercomplex variables in the seminal paper by Brackx, De Schepper, and Sommen in the mid-2000s (see \cite{BDS05}). It is noteworthy that both papers unjustly overlooked the interplay between Bernstein-type spaces and Bernstein's inequalities. This paper can be seen as a comprehensive exploration of these concepts.
Inspired by these works, we will also conduct a series of investigations to enable us to formulate the notion of Bernstein-type spaces in the context of hypercomplex analysis.

\subsection{Synopsis}

This paper addresses a hypercomplex version of the Paley-Wiener theorem within the framework developed by Li, McIntosh, and Qian \cite{li1994clifford}, and McIntosh \cite{mcintosh1996clifford}. For a thorough review, see the book by Qian and Li \cite{QianLi19}. Our approach also integrates the work of Andersen et al. \cite{Andersen04, Andersen06, AndersenJeu10, Andersen14} and Pesenson et al. \cite{Pesenson98, Pesenson00, Pesenson01, Pesenson08, PesensonZayed09, Pesenson15}. 

The structure of the paper is as follows. In Section \ref{sec2}, we revise some preliminaries on Clifford algebras and Clifford-valued function spaces.
In Section \ref{sec4}, we introduce the fundamental tools of our approach, including the concept of Dirac operators and generalized Hardy and Cauchy-type kernels. {We also explain our choice of the Riesz-Hilbert transform and, consequently, the Dirac operator (Subsection \ref{subsec41}). 
The two main components of our approach are the following: First, the introduction of the space-fractional Dirac operator 
\begin{equation*}%\tag{$\mathcal{D}_{\theta}^\alpha$}
	\D_\theta^\alpha=e^{-\frac{i \pi\theta}{2}}(-\Delta)^{\frac{\alpha}{2}}\left(\cos\left(\frac{ \pi\theta}{2}\right)+i\sin\left(\frac{\pi\theta}{2}\right)\mathcal{H}\right),
\end{equation*}
where $\mathcal{H}$ denotes the Riesz-Hilbert transform. This operator  includes not only the Dirac operator $\D$ (see \eqref{DiracOpPDO}), but also the Riesz derivative $(-\Delta)^{\frac{\alpha}{2}}$ (case $\theta=0$)
}
 of order $\alpha$ and skewness $\theta$, with $0<\alpha\leq 1$ and $\vert 1-\theta\vert<\frac{\alpha}{2}$. This extends the pseudo-differential representation of the Dirac operator {$\D$ (see \eqref{DiracOp} in Subsection \ref{subsec41})} as previously considered in the first author's paper \cite{Bernstein2016}; secondly, the formulation of generalized Hardy and Bernstein-type spaces derived from the Clifford-valued solutions of the subsequent Cauchy problem \eqref{CauchySpaceFractionalDirac} in $\BR^{n+1}$.
%:
%\begin{equation}\tag{$\mathcal{CP}_{\theta}^\alpha$}
%	\label{CauchySpaceFractionalDirac}\left\{\begin{array}{lll} 
%		\displaystyle 	\partial_{x_0}\u(\X,x_0)+\D_\theta^\alpha \u(\X,x_0)=0 &,~\X\in \BR^n&,~x_0\in \BR \setminus \{0\} 
%		\\ \ \\
%		\displaystyle \u(\X,0)=\f(\X)&,~\X \in \mathbb{R}^n \\ \ \\
%		\displaystyle	\lim_{\vert x_0 \vert \rightarrow +\infty}\u(\X,x_0)=0&,~\X \in \mathbb{R}^n,
%	\end{array}\right.
%\end{equation}
%where $\partial_{x_0}:=\frac{\partial}{\partial x_0}$ denotes the partial derivative on the $x_0-$direction.

Following the introduction of a new class of Dirac-type operators, the subsequent step entails the development of a comprehensive array of novel techniques that have not been addressed by the monogenic theory of Hardy spaces (see \cite{Mitrea94}). This necessity arises primarily from the absence of a Cauchy integral type formula in the space-fractional setting. 
To address this critical gap, the generalized Hardy spaces $H_{\alpha,\beta}^p\left(\BR_\pm^{n+1};\cl_{0,n}\right)$ and the underlying generalized Cauchy kernels $\E_{\alpha,n}^\pm$, have been defined and developed in Subsection \ref{subsec42} and Subsection \ref{subsec43}, respectively.

In developing such framework, we have integrated ideas from the following three contemporary research areas:
\begin{enumerate}
	\item[\textbf (I)] On the theory of {\textbf semigroups}, the fractional Laplacian of order $\alpha$, $-(-\Delta)^{\frac{\alpha}{2}}$, is the generator of the {\textit Feller semigroup} 
	$$ \left\{\exp\left(-te^{\frac{i\pi \gamma}{2}}(-\Delta)^{\frac{\alpha}{2}}\right)\right\}_{t\geq 0}, \text{ whereby  $0<\alpha\leq 1$ and $\vert \gamma\vert< \alpha$.}$$
%whereby $0<\alpha\leq 1$ and $\vert \gamma\vert< \alpha$.
	\item[\textbf (II)] On the theory of {\textbf probability}, the so-called {\textit L\'evy-Feller diffusion process} is uniquely determined, for values of $0<\alpha\leq 1$, by the kernel function
	\begin{equation*}%\tag{$\mathcal{K}_{\alpha,n}$}
	{\bm K}_{\alpha,n}(\X,z)=\frac{1}{(2\pi)^{n}}\int_{\BR^n}e^{-z\vert \XI \vert^\alpha} e^{i\langle \X,\XI\rangle}\dxi ,~~\X\in \BR^n,~~\Re(z)>0.
	\end{equation*}
	\item[\textbf (III)] In the theory of {\textbf singular integral operators}, the {\textit Riesz-Hilbert transform}  $\mathcal{H}$
	permits us to consider, for every $\f\in L^p(\BR^n;\cl_{0,n})$ ($1<p<\infty$), the functions $\displaystyle \f_\pm=\frac{1}{2}(\f\pm \mathcal{H}\f)$ as the boundary values of Hardy type spaces in the upper and lower half spaces, $\BR_+^{n+1}$ resp. $\BR_-^{n+1}$.
\end{enumerate}

We refer to \cite[Chapter IX]{Feller71} and \cite[Chapter 4]{Jacob01} for an overview of {\textbf (I)}, to \cite{BG1960} and the references therein for an overview of {\textbf (II)}, and to \cite{Mitrea94,mcintosh1996clifford,QianLi19} for an overview of {\textbf (III)}.
The combination of {\textbf (I)}, {\textbf (II)}, and {\textbf (III)} in this paper provides a comprehensive framework for dealing with function spaces and distributions. This framework is not limited to power series representations of solutions to the Cauchy problem \eqref{CauchySpaceFractionalDirac}, which several authors have previously explored using the Cauchy-Kovalevskaya extension (cf. \cite{Delanghe2001}). 

\begin{remark}
Although statement {\textbf (III)} holds in the borderline cases $p=1$ and $p=\infty$ (cf.~\cite{DangMaiQian2020}), we will focus our analysis in this paper on the range $1<p<\infty$ to avoid further technicalities related to arguments of weak boundedness in $BMO$ type spaces (see~\cite{DuongYan05} and the references therein for an overview). We plan to pursue these cases in a separate paper.
\end{remark}

In Section \ref{sec6}, we will consider the generalization of Paley-Wiener theory to the hypercomplex setting, focusing on the solutions of the equation \eqref{CauchySpaceFractionalDirac} that are generated from the one-parameter semigroup $\{e^{-x_0\D_\theta^\alpha}\}_{x_0\in \BR}$ of exponential type $R^\alpha$. 
The central tenet underpinning our construction is the investigation of the growth behavior of the sequences of functions $\left(\left(~\D_\theta^\alpha~\right)^k\f_\pm\right)_{k\in \BN_0}$, which play the role of the sequence of functions $\left(f^{(k)}\right)_{k\in \BN_0}$ generated by the derivatives of $f$ in the one-dimensional setting.
It is important to note that these two sequences of Clifford-valued functions are generated from the iterated powers of $\D_{\theta}^\alpha$ and the boundary values $\displaystyle \f_\pm$ of the Hardy spaces, as mentioned in {\textbf (III)}.  

 In real Paley-Wiener theory, one typically starts by choosing a real-valued bump function $\PSI$ such that $\widehat{\PSI}(\XI)=1$ in $\operatorname{supp}\widehat{\f}$, so that each $\f_\pm$ can be represented, in the distributional sense, as $\f_\pm=\PSI *\f_\pm$.
Such a fact is inherited from the framework considered by Andersen and Jeu in \cite{AndersenJeu10} and has been used to simplify parts of the proofs available in most of the literature, including the proofs on the seminal papers by Bang \cite{Bang90} and Tuan \cite{Tuan99}. In our concrete case, this serves as a starting point, first in Section \ref{subsec61}, Theorem \ref{RealPaleyWiener-BumpFunction}, to reformulate Andersen's proof \cite{Andersen06} on the Hankel transform for Clifford-valued functions which are {\textit radially symmetric}, and later, in Subsection \ref{subsec62}, to extend the Paley-Wiener framework presented by Pesenson-Zayed in \cite{PesensonZayed09} and by Pesenson in \cite{Pesenson15} to the hypercomplex setting. 

Although in the case of {\textit radially symmetric functions}, the normalized Bessel functions of order $\frac{n}{2}$, $\PSI(\X)=\left(R\vert \X\vert \right)^{-\frac{n}{2}}J_{\frac{n}{2}}\left(R\vert \X\vert\right)$, seem to be natural candidates for defining compactly supported functions on $\overline{B(0,R)}$ satisfying $\widehat{\PSI}(\XI)=1$ in $\operatorname{supp}\widehat{\f}$ (cf. ~\cite[p.~118]{Jacob01}), it should be emphasized that the existence of $\PSI$ is topologically guaranteed by the so-called Uryshon's Lemma (see \cite[Exercise 1.15]{LiebLoss2001}).

The main result of the paper is in Section \ref{subsec62}, Theorem \ref{PaleyWienerEquivalences}, where we prove a chain of Paley-Wiener equivalences, which in several papers are called real Paley-Wiener theorems.

As will be demonstrated in Section \ref{Applications}, the tools developed in Section \ref{sec6} can be used for two purposes: first, to generalize the proof given by Kou and Qian in their paper \cite{KouQian2002} and to deepen the characterization of Bernstein spaces as presented in \cite[Theorem 5.1]{FranklinHoganLarkin17} (see Theorem \ref{PW-KouQian-Generalization}); Secondly, the tools developed herein can provide more precise information about $\operatorname{supp}\widehat{\f}$ (see Theorem \ref{SpectralRadiusTheorem}).

In the final Section \ref{Conclusion}, we draw some conclusions and discuss open problems.

\begin{remark}
In the case where $p=2$, the compactness constraints imposed throughout the paper are not required a priori due to the direct application of Plancherel's theorem (see, for example, \cite[Theorem 3.5]{AndersenJeu10}). This observation is consistent with the proofs of the Paley-Wiener theorem obtained in the works \cite{KouQian2002,FranklinHoganLarkin17}, where the condition $\operatorname{supp}\widehat{\f}\subseteq \overline{B(0, R)}$ was replaced by $\operatorname{supp}\widehat{\f}\subseteq {B(0, R)}$.
\end{remark}

% In Section \ref{sec6} we prove a chain of Paley-Wiener equivalences, which in several papers are called real Paley-Wiener theorems. In Section \ref{Applications} we consider two applications of the main result obtained in Section \ref{sec6}: First, we generalize the Paley-Wiener theorems proved by Kou and Qian (see ~\cite{KouQian2002}) and Franklin, Hogan and Larkin (see ~\cite{FranklinHoganLarkin17}). We then explicitly compute the radius of the maximal ball $B(0,R)$ containing $\mbox{supp}\widehat{\f}$. In Section \ref{Conclusion} we conclude with an outlook on our approach.

\section{Preliminaries}\label{sec2}
\subsection{Clifford algebras}\label{subsec21}
Let $\cl_{0,n}$ be the universal Clifford algebra of signature $(0,n)$, generated by the basis $\e_1,\e_2,\ldots,\e_n$ of $\BR^n$. The multiplication in $\cl_{0,n}$ is induced by the graded anti-commuting relations
\begin{eqnarray}
	\label{CliffordBasis}\e_j\e_k + \e_k\e_j =-2\delta_{jk},  &\mbox{for any}& j,k=1,2,\ldots,n.
\end{eqnarray}

Now, let $\e_A=\e_{j_1}\e_{j_2}\ldots \e_{j_r}$ be a $r-$multivector of $\cl_{0,n}$ labeled by the ordered subset $A=\{j_1,j_2,\ldots,j_r\}$ of $[n]:=\{ 1,2,\ldots,n\}$, namely 
\begin{center}
	$1\leq j_1<j_2<\ldots<j_r\leq n$. 
\end{center}

In particular, vectors $\X=(x_1,x_2\ldots,x_n)$ and $\XI=(\xi_1,\xi_2,\ldots,\xi_n)$ of $\BR^n$ are represented through of the linear combinations
\begin{center}
 $\displaystyle \X=\sum_{j=1}^n x_j \e_j \mbox{ resp. }  \displaystyle \XI=\sum_{j=1}^n \xi_j \e_j,$ 
\end{center}
and the Euclidean inner product $\langle \X,\XI \rangle$ between $\X,\XI \in \BR^n$: 
\begin{center}
$\displaystyle \langle \X,\XI \rangle=\sum_{j=1}^n x_j\xi_j \in \BR $
\end{center}
can be casted in terms of anti-commutator identity, i.e. $\displaystyle \langle \X,\XI \rangle=-\frac{1}{2}(\X\XI+\XI \X)$.

By the preceding relation there holds that $\XI^2$ is a real number satisfying $\XI^2=-\vert \XI\vert^2$, where $\vert \XI\vert :={\langle \XI,\XI \rangle}^{\frac{1}{2}}$ stands for the Euclidean norm in $\BR^n$.
More generally, following \cite[Chapter 1]{BDS82}, one can endow an inner product by defining it for any ${\bm \lambda}, {\bm \mu} \in \cl_{0,n}$ by 
\begin{center}
$\displaystyle ({\bm \lambda} , {\bm \mu})_0 := 2^n [{\bm \lambda} {\bm \mu}^{\dagger}]_0 = 2^n \sum_A \lambda_A\mu_A , $
\end{center}
with $\displaystyle {\bm \mu}^{\dagger} = \sum_{A} \mu_A\overline{e_A}$, so that the corresponding norm is given by 
\begin{center}
$\displaystyle \vert{\bm \lambda}\vert_0 := \sqrt{({\bm \lambda} , {\bm \lambda})_0} = 2^{\frac{n}{2}} \sqrt{ [{\bm \lambda} {\bm \lambda}^{\dagger}]_0}$. %= 2^{\frac{n}{2}} \left(\sum_A \vert {\bm \lambda}_A\vert^2\right)^{\frac{1}{2}}.$
\end{center}

Hence, $\cl_{0,n}$ is simultaneously a real/complex Hilbert space and a Banach algebra. For the $\vert \cdot \vert_0-$norm, defined as above, one has the identity $\vert \e_0\vert_0 = 2^{\frac{n}{2}} \not = 1$. Additionally, the following lemma follows straightforwardly.
\\[0.1ex]
\begin{lemma}\label{CliffordIneqs}
For each ${\bm \lambda},{\bm \mu}\in \cl_{0,n}$, one has the following:
\begin{enumerate}
\item[(i)] 
({Submultiplicativity}) $
	\vert {\bm \lambda} {\bm \mu}\vert_0 \leq \vert {\bm \lambda}\vert_0 \vert{\bm \mu}\vert_0
$ \\[-2ex]
\item[(ii)] ({Triangle Inequality})
	$\vert{\bm \lambda}  + {\bm \mu}\vert_0  \leq \vert{\bm \lambda}\vert_0 + \vert{\bm \mu}\vert_0 .$ % \label{triangle}

\end{enumerate}
\end{lemma}
%For a sake of reader's convenience, the proof of Lemma \ref{CliffordIneqs} can be found in Appendix \ref{appendixA1}.

\begin{remark} 
\begin{itemize}
\item[]
\item[(i)]	The condition $\vert \e_0\vert_0 = 2^{\frac{n}{2}}$ is needed to get a Banach algebra. According to L. Ingelstam \cite{ingelstam}, a Hilbert space with an inner product $(\X,\Y)$ is also an associative Banach algebra with identity $\e_0$, and having a norm $\|\X\|=(\X,\X)^{\frac{1}{2}}$ satisfying $\|\X \Y\| \leq \|\X\| \cdot \|\Y\|$ and $\|\e_0\|=1$, is called a Hilbert algebra with identity $\e_0$. A real Hilbert algebra with identity is isomorphic to the real numbers, the complex numbers, or the quaternions. This result excludes the possibility that Clifford algebras, generally, can be Hilbert algebras with identity.
\item[(ii)] The standard absolute value of elements in the Clifford algebra $\displaystyle \vert \lambda \vert =\sum_A \vert \lambda_A\vert^2 $ doesn't fulfill the submultiplicativity condition. Instead of using the factor $2^{\frac{n}{2}}$ we could also define a Clifford operator norm \cite{GilbertMurray1991}. 
\item[(iii)] The submultiplicativity is needed to prove simple properties of function spaces such as H\"older's or Young's inequality.
\end{itemize}
\end{remark}

As a result, $\vert\cdot \vert_0$ has the properties of an absolute value. From now on, we will use this absolute value to define Clifford-valued function spaces.

\subsection{Function spaces}\label{sec3}
%\subsubsection{$L^p-$spaces}\label{subsec31}
\begin{definition}[\cite{BDS82}] $X_{(r)}$ is a \textbf{unitary right $\cl_{0,n} $-module}, when $(X_{(r)}, +)$ is an abelian group and the mapping $(\f, {\bm \lambda}) \to \f{\bm \lambda}$ from $ X_{(r)}\times\cl_{0,n}  \to X_{(r)}$ is defined such that for all ${\bm \lambda},{\bm \mu}\in \cl_{0,n}$ and $\f,\g\in X_{(r)}:$
	\begin{enumerate}
		\item $\f({\bm \lambda} + {\bm \mu}) = \f{\bm \lambda} + \f{\bm \mu},$
		\item $\f({\bm \lambda}{\bm \mu}) = (\f{\bm \lambda}){\bm \mu}, $
		\item $(\f+\g){\bm \lambda} = \f{\bm \lambda} + \g{\bm \lambda}, $
		\item $\f\e_0 = \f.$
	\end{enumerate}
\end{definition}

\begin{definition}\label{LpSpaces}
	A Clifford-valued function $\displaystyle \f(\X ) = \sum_A f_A(\X)\e_A$ belongs to \hfill \\ $L^p(\mathbb{R}^n; \cl_{0,n}),$ $ 1\leq p <\infty,$ if and only if 
	$$ \|\f\|_p : = \left(\int_{\mathbb{R}^n} \vert \f(\X)\vert^p_0 \dx\right)^{\frac{1}{p}} < \infty . $$
\end{definition}

\begin{remark}
\begin{itemize}\itemsep1ex
\item[~]
\item[(i)]
Obviously, all $L^p-$spaces are  unitary right $\cl_{0,n}-$modules.
	It should be mentioned that $\| \cdot \|_p$ is a real norm. But it is only a Clifford pseudonorm, because \hfill \\ $\| {\bm \lambda}{\f}\|_p \leq \vert {\bm \lambda} \vert_0 \| {\f} \|_p$ for ${\bm \lambda}\in \cl_{0,n}.$
\item[(ii)] The norms $\|\f\|_p$ and $\displaystyle \left(\sum_A \int_{\mathbb{R}^n} \vert f_{A}(\X)\vert^p \dx \right)^{\frac{1}{p}}$ are equivalent, i.e.,
		a Clifford-valued function $\f$ belongs to $L^p(\mathbb{R}^n; \cl_{0,n}),$ \\ $1\leq p<\infty,$ if and only if each component $f_A$ belongs to $L^p(\mathbb{R}^n).$
\item[(iii)] Functions of these Clifford analogues to the $L^p$-spaces fulfil the usual Minkowski's, H\"older's and Young's inequality.
\end{itemize}		
\end{remark}

%This norm fulfills the usual inequalities:
%
%\begin{proposition}\label{LpNormProp}
%	\begin{itemize}\itemsep1ex
%		\item[~]
%		\item[(i)] (\textbf{Equivalent norms}) The norms $\|\f\|_p$ and $\displaystyle \left(\sum_A \int_{\mathbb{R}^n} \vert f_{A}(\X)\vert^p \dx \right)^{\frac{1}{p}}$ are equivalent, i.e.,
%		a Clifford-valued function $\f$ belongs to $L^p(\mathbb{R}^n; \cl_{0,n}),$ \\ $1\leq p<\infty,$ if and only if each component $f_A$ belongs to $L^p(\mathbb{R}^n).$
%		\item[(ii)](\textbf{Minkowski's inequality})\label{MinkowskiIneq}
%		For $1<p<\infty$ and $\f,\g\in L^p(\mathbb{R}^n; \cl_{0,n})$
%		$$
%		\| \f + \g \|_p \leq \| \f \|_p+\| \g \|_p. 
%		$$
%		\item[(iii)](\textbf{H\"older's inequality})\label{HolderIneq}
%		Let $p,q > 0$ such that $\frac{1}{p} + \frac{1}{q} = 1$ and  $\f\in L^p(\mathbb{R}^n; \cl_{0,n})$, $\g\in L^q(\mathbb{R}^n; \cl_{0,n})$, then
%		$$ \| \f \g \|_1 =   \int_{\mathbb{R}^n} \vert f(\X)\g(\X)\vert_0 \dx \leq \| \f \|_p \| \g \|_q . $$
%		\item[(iv)](\textbf{Young's inequality})\label{YoungIneq} Let $p,q,r \in \mathbb{R}$, $p,q,r\geq 1,$ and 
%		$$ \frac{1}{p} + \frac{1}{q}=1+\frac{1}{r}. $$
%		Let $\f\in L^p(\mathbb{R}^n; \cl_{0,n}), \g\in L^q(\mathbb{R}^n; \cl_{0,n})$, then the convolution \\ $\f*\g \in L^r(\mathbb{R}^n; \cl_{0,n})$
%		and the following inequality is satisfied
%		$$ \|\f*\g\|_r \leq  \|\f\|_p\,\| \g\|_q . $$
%	\end{itemize}
%\end{proposition}

Next, we turn our attention to the inner product structure on a unitary right $\cl_{0,n}$-module.
\begin{definition}[\cite{BDS82}, p. 13]
	Let $H_{(r)}$ be a unitary right $\cl_{0,n}-$module. Then, the function 
	$\langle\cdot, \cdot \rangle: H_{(r)}\times H_{(r)} \to \cl_{0,n} $ defines an \textbf{inner product} on $H_{(r)}$ if for all $\f,\g,{\bm h} \in H_{(r)}$ and ${\bm \lambda} \in \cl_{0,n},$ 
	\begin{enumerate}
		\item $\langle \f, \g+{\bm h}\rangle = \langle \f,\g\rangle + \langle \f,{\bm h}\rangle$
		\item $\langle \f,\g{\bm \lambda}\rangle = \langle \f,\g\rangle {\bm \lambda} $
		\item $\langle \f,\g\rangle  = \langle \g,\f\rangle^{\dagger}$
		\item $\langle \f,\f\rangle_0\in\mathbb{R}_0^+$ and $\langle \f,\f\rangle _0=0$ if and only if $\f=0$
		\item $\langle \f{\bm \lambda},\f{\bm \lambda}\rangle _0\leq\vert{\bm \lambda}\vert_0^2\langle \f,\f\rangle _0.$
	\end{enumerate}
	The accompanying \textbf{norm} on $H_{(r)}$ is $\left\|\f\right\|^2=\langle \f,\f\rangle_0$.
\end{definition}

Hence $L^2(\mathbb{R}^n; \cl_{0,n})$ is a unitary right $\cl_{0,n}$-module with inner product
$$ \langle \f, \g\rangle:= 2^n \int_{\mathbb{R}^n} \f(\X)^{\dagger} \g(\X) \dx,$$
 and the induced norm is 
$$ \|\f\|^2 :=  \langle \f, \f\rangle_0 = 2^n \int_{\mathbb{R}^n} [\f(\X)^{\dagger} \f(\X)]_0 \dx = \int_{\mathbb{R}^n} \vert \f(\X)\vert^2_0 \dx = \|\f \|_2^2 $$
%{Due to the Clifford structure of our spaces, we obtain a slightly modified Cauchy-Schwarz inequality:}
%\begin{proposition}[Cauchy-Schwarz inequality]\label{CauchySchwarzIneq}
%	If $\langle\cdot,\cdot\rangle$ is an inner product on a unitary right $\cl_{0,n}$-module $H_{(r)}$ and $\left\|\f\right\|^2=2^n\langle \f,\f\rangle_0$ then
%	\[\vert \langle \f,\g\rangle_0\vert\leq 2^{n}\left\|\f\right\|\,\left\|\g\right\|\]
%	holds for all $\f,\g\in H_{(r)}$.
%\end{proposition}
%
%
%For this specific inner product, we obtain the classical Cauchy-Schwarz inequality:
%\[\vert\langle \f,\g\rangle\vert_0\leq \left\|\f\right\|_2\,\left\|\g\right\|_2\]
%for all $\f,\g\in L^2(\mathbb{R}^n; \cl_{0,n})$.

%\subsubsection{Schwartz spaces}\label{subsec32}
We start to define the Schwartz spaces as in \cite{grochenig}.  Given a multi-index $ \alpha = (\alpha_1, \ldots, \alpha_n) \in \mathbb{N}_0^n,$ we take the standard multi-index notations
\begin{center}
 $\displaystyle \vert \alpha \vert = \sum_{j=1}^n \alpha _j$,   $\displaystyle \X^{\alpha}=\prod_{j=1}^d x_j^{\alpha_j}$, 
\end{center}
and the shorthand notations 
\begin{eqnarray}
\displaystyle \partial^{\alpha} = \frac{\partial^{\alpha_1}}{\partial x_1^{\alpha_1}}\cdots  \frac{\partial^{\alpha_n}}{\partial x_n^{\alpha_n}} &\text{and}& \displaystyle X^\alpha \f(\X ) = \X^\alpha \f(\X)
\end{eqnarray}
 for multi-index partial derivative and multiplication operator of order $\vert \alpha \vert$, respectively.  
% The Schwartz space  $\mathcal{S}(\BR^n)$ consists of all real-valued $C^{\infty}$-functions $f$ on $\BR^n$ such that 
%$$ \sup_{\X\in\BR^n} \vert \partial^{\alpha}X^{\beta}f(\X)\vert < \infty \quad \text{for all}\quad \alpha, \beta \in \mathbb{N}_0.$$
%Therefore, it is natural to 
We define the Clifford Schwartz-type space  $\mathcal{S}(\BR^n;\cl_{0,n})$ as the set of all Clifford-valued $C^{\infty}$-functions $\f$ on $\BR^n$ such that 
\begin{eqnarray*}
 \sup_{\X\in\BR^n} \vert \partial^{\alpha}X^{\beta}\f(\X)\vert_0 < \infty \quad \text{for all}\quad \alpha, \beta \in \mathbb{N}_0. \\[0.5ex]
\end{eqnarray*}

\begin{proposition}\label{SchwartzSpaceProposition}
The function $\f$ defined by $\displaystyle \f(\X)= \sum_A f_A(\X) e_A$ belongs to $\mathcal{S}(\BR^n;\cl_{0,n})$ if and only if all $f_A \in \mathcal{S}(\BR^n).$ 
\end{proposition} 

%\textcolor{orange}{\sout{See Appendix \ref{appendixA4} for the proof of Proposition \ref{SchwartzSpaceProposition}.}}

\section{Clifford analysis}\label{sec4}

\subsection{Dirac operators of Riesz-Feller type}\label{subsec41}

The goal of this subsection is to introduce a hypercomplex analog of the so-called Riesz-Feller derivative, which corresponds to the space-fractional analog of the Dirac operator \\[-5ex]
\begin{eqnarray}
	\label{DiracOp}
	\D=\sum_{j=1}^n\e_j\partial_{x_j}.
\end{eqnarray}
To approach this properly, let us begin by defining the Fourier transform of $\f\in \mathcal{S}(\BR^n;\cl_{0,n})$ and the Fourier inverse of $\g\in \mathcal{S}(\BR^n;\cl_{0,n})$ by
\begin{eqnarray}
	%\label{FourierTransform}
	\widehat{\f}(\XI):=(\mathcal{F} \f)(\XI)=\int_{\BR^n} \f(\X)e^{-i \langle \X,\XI \rangle}\dx, & \XI \in \BR^n, \\[1ex]
%\end{eqnarray}
% and by
%\begin{eqnarray}
	%\label{FourierInverse}
(\mathcal{F}^{-1} \g)(\X)=\frac{1}{(2\pi)^n}\int_{\BR^n} \g(\XI)e^{i \langle \X,\XI \rangle}\dxi, & \X \in \BR^n.
\end{eqnarray}
%the Fourier inverse of $\g\in \mathcal{S}(\BR^n;\cl_{0,n})$.

The mapping property $\mathcal{F}:\mathcal{S}(\BR^n;\cl_{0,n})\longrightarrow \mathcal{S}(\BR^n;\cl_{0,n})$ can be extended to $L^2(\BR^n;\cl_{0,n})$ by using density arguments, and to the {\textit space of tempered distributions} $\mathcal{S}'(\BR^n;\cl_{0,n})$ by duality arguments.
%, via the duality relation
%\begin{eqnarray*}
%	\langle \widehat{\f},\g \rangle=\langle \f,\widehat{\g}\rangle, & \f\in \mathcal{S}(\BR^n;\cl_{0,n}), & \g\in \mathcal{S}'(\BR^n;\cl_{0,n}).
%\end{eqnarray*}

To draw a parallel with the approach of the first author \cite{Bernstein2016}, let us closely examine the Dirac operator $\D$, defined by \eqref{DiracOp}, and the Laplace operator 
\begin{center}
$\displaystyle \Delta=\sum_{j=1}^{n}\partial_{x_j}^2.$
\end{center}

First, we remember that $\D$ and $\Delta$ are connected by the factorization property $\D^2 = -\Delta$. Moreover, we can extend the actions of both operators to spaces of functions and distributions. Namely, for $\f\in \mathcal{S}(\BR^n;\cl_{0,n})$ holds
\begin{eqnarray}
	\label{spectralDLaplace}
	%\begin{array}{lll}
		\mathcal{F}(\D \f)(\XI) =-i\XI \widehat{\f}(\XI) \text{ and }
		\mathcal{F}(-\Delta\f)(\XI) ={\vert \XI\vert^2} \widehat{\f},
		(\XI),&  \XI\in \BR^n.
	%\end{array}
\end{eqnarray}

The spherical decomposition of the Fourier symbol of $\D$:
\begin{eqnarray*}
	-i\XI=\vert \XI\vert \frac{-i\XI}{\vert\XI\vert} &,&\XI \in \BR^n\setminus \{\bm 0\},
\end{eqnarray*}
together with the set of identities $(-i\XI)^2=\vert \XI\vert^2$ and $\vert \XI\vert^\alpha=\left(~\vert \XI\vert^2~\right)^\frac{\alpha}{2}$,
allows us to determine the pseudo-differential operators $(-\Delta)^{\frac{\alpha}{2}}$ and $\mathcal{H}:=\D(-\Delta)^{-\frac{1}{2}}$ in terms of the Fourier symbols $\displaystyle \vert \XI\vert^\alpha$ and $\displaystyle \frac{-i\XI}{\vert \XI\vert}$, respectively. More precisely, for $\f\in \mathcal{S}(\BR^n;\cl_{0,n})$ the component-wise actions of $(-\Delta)^{\frac{\alpha}{2}}\f$ and $\mathcal{H}\f$ are uniquely determined by the spectral formulas
\begin{eqnarray}
	\label{RieszDerivative}	 \mathcal{F}\left(\left(-\Delta\right)^{\frac{\alpha}{2}}\f\right)(\XI) ={\vert \XI\vert^\alpha} \widehat{\f}
	(\XI), &\XI\in \BR^n, \\
	\label{RieszHilbert}		\mathcal{F}(\mathcal{H} \f)(\XI) =\frac{-i\XI}{\vert \XI\vert} \widehat{\f}(\XI), &\XI\in \BR^n \setminus \{0\}. 
\end{eqnarray}

That is equivalent to $(-\Delta)^{\frac{\alpha}{2}}:=\mathcal{F}^{-1}{\vert \XI\vert^\alpha}\mathcal{F}$ resp. $\displaystyle \mathcal{H}=\mathcal{F}^{-1}\frac{-i\XI}{\vert \XI\vert}\mathcal{F}$.

The operator $(-\Delta)^{\frac{\alpha}{2}}$ corresponds to the Riesz derivative of order $\alpha$, while $\mathcal{H}$ represents the hypercomplex generalization of the one-dimensional Riesz transform, due to the unitary property ${\bm \omega}^2=1$ of $\displaystyle {\bm \omega}:=\frac{-i\XI}{\vert \XI\vert}$. We will refer to $\mathcal{H}$ as the {\textit Riesz-Hilbert transform} (cf. \cite{Bernstein2017}, \cite{HRS-2014}).

Based on this, the Dirac operator \eqref{DiracOp} allows the following pseudo-differential reformulation\\[-6ex]
\begin{eqnarray}
	\label{DiracOpPDO}
	\D=\vert \D \vert \operatorname{sgn}(\D),
\end{eqnarray}
where $\vert \D\vert:=(-\Delta)^{\frac{1}{2}}$ is the radial part of $\D$ and $\mbox{sgn}(\D):=\mathcal{H}$ is the phase part of $\D$ (cf.~\cite{mcintosh1996clifford}).

Using the Fourier convolution theorem, both operators can also be expressed as singular integral operators.  For $\mathcal{H}$, it is important to emphasize the closed formula (see \cite{li1994clifford, mcintosh1996clifford}):
\begin{eqnarray} \label{RieszHilbertSingularInt} \mathcal{H} \f(\X)=\frac{\Gamma\left(\frac{n+1}{2}\right)}{\pi^{\frac{n+1}{2}}} P.V. \int_{\BR^n}\frac{\X-\Y}{\vert \X-\Y\vert^{n+1}}\f(\Y)\dy, &\X \in \BR^n \end{eqnarray}

This follows from the fact that $\mathcal{H}$ is given by the linear combination
\begin{eqnarray} \label{RieszHilbertLinearCombination} \mathcal{H}=\sum_{j=1}^{n}\e_j R_j,   \text{ where each } R_j:=\mathcal{F}^{-1}\frac{-i\xi_j}{\vert \XI\vert} \mathcal{F} \end{eqnarray} (a directional Riesz transform along the $x_j$ axis) admits a singular integral representation in terms of the kernel function (cf.~\cite[p.~224]{SteinWeiss72}):
\begin{eqnarray}
	\label{RieszKernel}
	E_j(\X)=\frac{\Gamma\left(\frac{n+1}{2}\right)}{\pi^{\frac{n+1}{2}}}\frac{x_j}{\vert\X \vert^{n+1}}&,~\X \in \BR^{n} \setminus \{\bm 0\} &  (j=1,2,\ldots,n).
\end{eqnarray}

%%% @nelson
%In addition, the operator $\D_\theta^\alpha$, which is defined by \eqref{DiracRieszFellerFourier}, allows the pseudo-differential representation
%\begin{equation}%\tag{$\mathcal{D}_{\theta}^\alpha$}
%	\label{RieszFellerOp}
%	\D_\theta^\alpha=e^{-\frac{i \pi\theta}{2}}(-\Delta)^{\frac{\alpha}{2}}\left(\cos\left(\frac{ \pi\theta}{2}\right)+i\sin\left(\frac{\pi\theta}{2}\right)\mathcal{H}\right),
%\end{equation}
%which includes not only the Dirac operator $\D$ (see \eqref{DiracOpPDO}), but also the Riesz derivative $(-\Delta)^{\frac{\alpha}{2}}$ (case $\theta=0$) and the fractional Riesz-Hilbert transform (limit case $\alpha\rightarrow 0^+$), studied in detail by the first author in \cite{Bernstein2016,Bernstein2017}.

\begin{definition}
For functions $\f \in \mathcal{S}(\BR^n;\cl_{0,n})$, the so-called {\textit Dirac operator of Riesz-Feller type} is defined by the spectral formula
\begin{eqnarray}
	\label{DiracRieszFellerFourier}
	\mathcal{F}\left(\D_\theta^\alpha \f\right)(\XI)=\vert \XI \vert^\alpha h_{\theta}(\XI) \widehat{\f}(\XI), & \XI \in \BR^n\setminus \{\bm 0\}, \\[2ex]
%\end{eqnarray}
%\begin{eqnarray}
	\text{where  }	h_{\theta}(\XI):=e^{-\frac{i\pi \theta}{2}}\left(\cos\left(\dfrac{\pi \theta}{2}\right)+\dfrac{\XI}{\vert \XI\vert}\sin\left(\dfrac{\pi \theta}{2}\right)\right), & \XI \in \BR^n\setminus \{\bm 0\}. \label{RieszHilbert-h}
\end{eqnarray}
Following the standard conventions on Riesz-Feller derivatives, we denote $\alpha$ as the order and $\theta$ as the skewness of $\D_\theta^\alpha$, such that
\begin{center}
	$0<\alpha\leq 1$~\&~$\displaystyle \vert 1-\theta\vert <\frac{\alpha}{2}$.
\end{center}
\end{definition}
In addition, the operator $\D_\theta^\alpha$, which is defined by \eqref{DiracRieszFellerFourier}, allows the pseudo-differential representation
\begin{equation}%\tag{$\mathcal{D}_{\theta}^\alpha$}
	\label{RieszFellerOp}
	\D_\theta^\alpha=e^{-\frac{i \pi\theta}{2}}(-\Delta)^{\frac{\alpha}{2}}\left(\cos\left(\frac{ \pi\theta}{2}\right)+i\sin\left(\frac{\pi\theta}{2}\right)\mathcal{H}\right),
\end{equation}
which includes not only the Dirac operator $\D$ (see \eqref{DiracOpPDO}), but also the Riesz derivative $(-\Delta)^{\frac{\alpha}{2}}$ (case $\theta=0$) and the fractional Riesz-Hilbert transform (limit case $\alpha\rightarrow 0^+$), studied in detail by the first author in \cite{Bernstein2016,Bernstein2017}.

The operator $\D_\theta^\alpha$ corresponds, up to the phase factor $e^{-\frac{i\pi\theta}{2}}$, to the space-fractional Dirac operator, already considered by the second author in \cite{Faustino21}. Throughout the paper, $\D_\theta^\alpha$ will be referred to as the \textit{Dirac operator of Riesz-Feller type} because of its several similarities to the one-dimensional Riesz-Feller derivative. For further details, we refer to the paper \cite{Mainardi01} by Mainardi-Luchko-Pagnini and the references therein.

The study of these operators and their variants, as infinitesimal generators of semigroups, goes back to the operator-theoretical approaches of Bochner \cite{Bochner49} and Feller \cite{Feller62}, connecting stable distributions with diffusion equations, and has been a source of many fundamental developments in the theory of Hardy spaces in a broad sense (cf.~\cite{Yan08}).
In addition to its similarity to the polar decomposition underlying complex numbers, it is noteworthy that the choice of $h_{\theta}(\XI)$ as the phase component of the Clifford vector ${\bm z}=\vert \XI \vert^\alpha h_{\theta}(\XI)$ holds particular significance in the field of optics. 
%%% @nelson
%The fractional Riesz-Hilbert transform we will use is constructed from the Riesz transform in the same way that this fractional Hilbert transform is obtained from the Hilbert transform, and it shares the same properties  \cite{Bernstein2017, HRS-2014, Bernstein-2018}. %\\

The definition of $h_{\theta}(\XI)$ by \eqref{RieszHilbert-h} is derived from the Fourier symbol of the Riesz-Hilbert transform $\mathcal{H}$. This approach parallels the methodology used to derive the fractional Hilbert transform from the Hilbert transform (see \cite{LTR-1997, LMZ-1998}). Furthermore, it extends the Fourier symbol of the fractional Hilbert transform in a coherent way, preserving its intrinsic properties \cite{Bernstein2017, HRS-2014, Bernstein-2018}. This makes it particularly valuable for analyzing the local orientation, amplitude, and phase of wave functions, as it seamlessly integrates the influences of both near and far fields.
The above aspects distinguish our work from much of the existing literature on fractional Dirac operators (see, e.g., \cite{Raspini01, MuslihAgrawalBaleanu10, QuanUhlmann24}). In our study, the optical implementation of the fractional Hilbert transform was crucial and served as the fundamental basis for the pseudo-differential representation of $\D_\theta^\alpha$, provided by \eqref{RieszFellerOp}. 
	We also emphasize that our methodology is in contrast to that proposed by Zayed in \cite{Zayed-1997}, where the fractional Fourier transform serves as the basis for constructing a fractional Hilbert transform. Furthermore, this technique is extended to the Riesz transform, allowing the derivation of fractional Riesz transforms; see, for example, \cite{FGLWY23}.

%%% @nelson
%This approach differs from that of Zayed in \cite{Zayed-1997}, where the fractional Fourier transform is used to construct a fractional Hilbert transform. This approach is also extended to the Riesz transform to obtain fractional Riesz transforms; see, for example, \cite{FGLWY23}.

\subsection{Generalized Hardy spaces $H^p_{\alpha,\beta}$}\label{subsec42}

Let us now consider the Cauchy problem as follows:
\begin{equation}\tag{$\mathcal{CP}_{\theta}^\alpha$}
	\label{CauchySpaceFractionalDirac}\left\{\begin{array}{lll} 
		\displaystyle 	\partial_{x_0}\u(\X,x_0)+\D_\theta^\alpha \u(\X,x_0)=0 &,~\X\in \BR^n&,~x_0\in \BR \setminus \{0\} 
		\\[1ex]
		\displaystyle \u(\X,0)=\f(\X)&,~\X \in \mathbb{R}^n, \\[1ex]
		\displaystyle	\lim_{\vert x_0 \vert \rightarrow +\infty}\u(\X,x_0)=0&,~\X \in \mathbb{R}^n,
	\end{array}\right.
\end{equation}
%where $\partial_{x_0}:=\frac{\partial}{\partial x_0}$ denotes the partial derivative on the $x_0-$direction.
%\eqref{CauchySpaceFractionalDirac} 

To study it from a Hardy space point of view, we need to define appropriate function spaces that capture the order and skewness parameters, $\alpha$ and $\theta$, of $\D_\theta^\alpha$ on the upper and lower half-spaces, $\BR_{+}^{n+1}$ and $\BR_{-}^{n+1}$, of $\BR^{n+1}$.
First and foremost, set 
\begin{eqnarray}
	\label{HardySpacesRn}
	L^{p,\pm}\left(\mathbb{R}^n;\cl_{0,n}\right) = \left\{\frac{1}{2}\left(\f\pm \mathcal{H}\f\right)~:~\f\in L^p(\BR^n;\cl_{0,n})\right\}.
\end{eqnarray}

We recall that 
\begin{eqnarray}
	\label{ProjectionOpHilbert}
	%\begin{array}{lll}
		\frac{1}{2}(I+\mathcal{H}) = \mathcal{F}^{-1} \chi_-(\XI) \mathcal{F}, \quad \frac{1}{2}(I-\mathcal{H}) = \mathcal{F}^{-1} \chi_+(\XI)  \mathcal{F},
	%\end{array}
\end{eqnarray}
where $\chi_-(\XI)$ and $\chi_+(\XI)$ represent the Fourier symbols
\begin{eqnarray}
	\label{ProjectionOp}
	\chi_\pm(\XI)=\frac{1}{2}\left(1\pm  \frac{i\XI}{\vert \XI\vert}\right), & \XI\in \BR^n\setminus \{\bm 0\}.
\end{eqnarray} 

Based on this, the following reformulation of \eqref{RieszFellerOp} in terms of $\frac{1}{2}\left(I\pm \mathcal{H}\right)$, given by
\begin{eqnarray}
	\label{RieszFellerHardy}
	\D_\theta^\alpha&=&(-\Laplace)^{\frac{\alpha}{2}}\frac{1}{2}(I+\mathcal{H})+e^{-i{\pi\theta}}(-\Laplace)^{\frac{\alpha}{2}} \frac{1}{2}(I-\mathcal{H})
\end{eqnarray}
is obtained by reformulating the Fourier symbol \eqref{RieszHilbert-h} as a linear combination of the symbols \eqref{ProjectionOp}. That is,
$h_\theta(\XI)=\chi_-(\XI)+e^{-i{\theta\pi}}\chi_+(\XI)$.

On the other hand, considering the properties (cf.~\cite[p. 671]{li1994clifford} \& \cite[subsection 5.2.1]{mcintosh1996clifford}) associated with the Fourier symbols 
\begin{eqnarray}
	\label{ProjectionOpProperties}
	\begin{array}{lll}
		\chi_+(\XI)+\chi_-(\XI)=1, &  \\
		\left(\chi_+(\XI)\right)^2=\chi_+(\XI), & \left(\chi_-(\XI)
		\right)^2=\chi_-(\XI),  	\\
		\chi_+(\XI)\chi_-(\XI)=0, & \chi_-(\XI)\chi_+(\XI)=0,
	\end{array}
\end{eqnarray}
it follows that $\displaystyle \frac{1}{2}(I\pm\mathcal{H})$ are projection operators, leading to the direct sum decomposition
$$ L^p\left(\mathbb{R}^n;\cl_{0,n}\right) =  L^{p,+}\left(\mathbb{R}^n;\cl_{0,n}\right)  \oplus  L^{p,-}\left(\mathbb{R}^n;\cl_{0,n}\right).$$
%satisfy the mapping property $\frac{1}{2}(I\pm\mathcal{H}):L^p(\BR^n;\cl_{0,n})\rightarrow L^p(\BR^n;\cl_{0,n})$, for values of $1<p<\infty$. 

This allows us to consider the $L^{p,\pm}-$spaces \eqref{HardySpacesRn} as hypercomplex analogues of the \textit{Hardy spaces} over $\BR^n$. Furthermore, with the help of the following analytic semigroups
\begin{eqnarray}
	\label{AnalyticSemigroups}	\left\{\exp\left(-x_0e^{i\pi\beta}(-\Delta)^{\frac{\alpha}{2}}\right)\right\}_{x_0\geq 0} &\mbox{resp.}&\left\{\exp\left(-x_0e^{-i\pi\beta}(-\Delta)^{\frac{\alpha}{2}}\right)\right\}_{x_0\leq 0}
\end{eqnarray}
we can extend it to the upper and lower half-spaces, $\BR_{+}^{n+1}$ and $\BR_{-}^{n+1}$, respectively, of $\BR^{n+1}$. This leads to the following definition:
\\[-0.5ex]
\begin{definition}[{\textbf Generalized Hardy spaces $H^p_{\alpha,\beta}$}]\label{Hardy-p}
	Let $0<\alpha\leq 1$, $\beta\in \BR$ and $1<p<\infty$ be given. We define $H^p_{\alpha,\beta}\left(\BR^{n+1}_+;\cl_{0,n}\right)$ and $H^p_{\alpha,\beta}\left(\BR^{n+1}_-;\cl_{0,n}\right)$
	as follows:
	\begin{enumerate}
		\item In case of $\vert \beta\vert < \dfrac{\alpha}{2}$, $\u_+\in H^p_{\alpha,\beta}\left(\BR^{n+1}_+;\cl_{0,n}\right)$ if, and only if:
		\begin{enumerate}
			\item[~]
			\item $\u_+$ is defined by
			\begin{eqnarray}
				\label{uPlus}
				\u_+(\X,x_0)=\exp\left({-x_0e^{i\pi \beta}(-\Delta)^{\frac{\alpha}{2}}}\right)\f_+(\X), &(\X,x_0)\in \BR^{n+1}_{+},
			\end{eqnarray} for some $\f_+\in L^{p,+}\left(\BR^n;\cl_{0,n}\right)$.
			\item[~]
			\item $\displaystyle \sup_{x_0>0} \|~\u_+(\cdot,x_0)~\|_p<\infty$.
		\end{enumerate}
		\item[~]
		\item In case of $\vert 1-\beta\vert < \dfrac{\alpha}{2}$, $\u_-\in H^p_{\alpha,\beta}\left(\BR^{n+1}_-;\cl_{0,n}\right)$ if, and only if:
		\begin{enumerate}
			\item[~]
			\item $\u_-$ is defined by
			\begin{eqnarray}
				\label{uMinus}
				\u_-(\X,x_0)=\exp\left({-x_0e^{-i\pi \beta}(-\Delta)^{\frac{\alpha}{2}}}\right)\f_-(\X), &(\X,x_0)\in \BR^{n+1}_{-},
			\end{eqnarray} for some $\f_-\in L^{p,-}\left(\BR^n;\cl_{0,n}\right)$.
			\item[~]
			\item $\displaystyle \sup_{x_0<0} \|~\u_-(\cdot,x_0)~\|_p<\infty$.
		\end{enumerate}
	\end{enumerate}
\end{definition}

\begin{remark}
	From the limit conditions
	\begin{eqnarray*}
		\lim_{x_0\rightarrow 0^+}	\u_+\left(\X,x_0\right)=\f_+(\X)&\mbox{\&}&\displaystyle \lim_{x_0\rightarrow 0^-}	\u_-\left(\X,x_0\right)=\f_-(\X)
	\end{eqnarray*}
	it holds that $\f_+=\frac{1}{2}\left(\f+ \mathcal{H}\f\right)$ and $\f_-=\frac{1}{2}\left(\f- \mathcal{H}\f\right)$ are the boundary values of the generalized Hardy spaces $H^p_{\alpha,\beta}\left(\BR^{n+1}_+;\cl_{0,n}\right)$ and $H^p_{\alpha,\beta}\left(\BR^{n+1}_-;\cl_{0,n}\right)$, respectively.
\end{remark}

The generalized Hardy spaces introduced above have a crucial property: they allow a complete characterization of the solutions of the Cauchy problem \eqref{CauchySpaceFractionalDirac}, as detailed below.
%%% @nelson
%A key property of the previously introduced generalized Hardy spaces is the characterization of the solutions to the Cauchy problem \eqref{CauchySpaceFractionalDirac}. 
\begin{proposition}\label{CauchyProblemHardy}
	For values of $0<\alpha\leq 1$, $\displaystyle \vert 1-\theta\vert< \frac{\alpha}{2}$ and $1<p<\infty$, we assume that $\f\in L^p(\BR^n;\cl_{0,n})$. We also assume that
	\begin{eqnarray*}
		\u_+\in H^p_{\alpha,0}\left(\BR^{n+1}_+;\cl_{0,n}\right) &\mbox{\&}&\u_-\in H^p_{\alpha,\theta}\left(\BR^{n+1}_-;\cl_{0,n}\right),
	\end{eqnarray*}
	where $H^p_{\alpha,0}\left(\BR^{n+1}_+;\cl_{0,n}\right)$ and $H^p_{\alpha,\theta}\left(\BR^{n+1}_-;\cl_{0,n}\right)$ correspond to the {\textit generalized Hardy spaces} introduced in Definition \ref{Hardy-p}.
	Then, the function $\u$ defined by
	\begin{eqnarray}
	\label{uAnsatz}
	\u(\X,x_0)=\begin{cases}
		\u_+(\X,x_0)&,~(\X,x_0)\in \BR^{n+1}_+ \\[1ex]
		\u_-(\X,x_0)&,~(\X,x_0)\in \BR^{n+1}_- \\[1ex]
		\displaystyle \lim_{x_0\rightarrow 0^+}\u_+(\X,x_0)+\lim_{x_0\rightarrow 0^-}\u_-(\X,x_0)&,~(\X,x_0)\in \BR^{n} \times \{0\}
	\end{cases}
\end{eqnarray}
	solves the Cauchy problem \eqref{CauchySpaceFractionalDirac} if, and only if, 
	$\u_{+}$ (case of $\beta=0$) and $\u_-$ (case of $\beta=\theta$)
are solutions of
	\begin{equation}%\tag{$CP_\alpha^+$}
		\label{CauchySpaceFractionalHplus}\left\{\begin{array}{lll} 
			\displaystyle 	\partial_{x_0}\u_{\pm}(\X,x_0)+e^{\pm i\pi\beta}(-\Delta)^{\frac{\alpha}{2}} \u_{\pm}(\X,x_0)=0 &, & (\X,x_0)\in \mathbb{R}^{n+1}_+
			\\[1ex]
			\displaystyle	\lim_{x_0 \rightarrow  0^{\pm}}\u_{\pm}(\X,x_0)=\tfrac{1}{2}\left(\f(\X) \pm \mathcal{H}\f(\X)\right)& , & \X \in \mathbb{R}^n \\[1ex]
			\displaystyle	\lim_{x_0 \rightarrow  +\infty}\u_{\pm}(\X,x_0)=0& , & \X \in \mathbb{R}^n
		\end{array}\right.
	\end{equation} 

%and
%	\begin{equation}%\tag{$CP_\alpha^-$}
%		\label{CauchySpaceFractionalHminus}\left\{\begin{array}{lll} 
%			\displaystyle 	\partial_{x_0}\u_-\left(\X,x_0\right)+e^{-i\pi\theta}(-\Delta)^{\frac{\alpha}{2}}\u_-\left(\X,x_0\right)=0 &, & (\X,x_0)\in \mathbb{R}^{n+1}_-
%			\\ \ \\
%			\displaystyle	\lim_{x_0 \rightarrow  0^-}\u_-\left(\X,x_0\right)=\frac{1}{2}\left(\f(\X)-\mathcal{H}\f(\X)\right)& , & \X \in \mathbb{R}^n \\ \ \\
%			\displaystyle	\lim_{x_0 \rightarrow  -\infty}\u_-\left(\X,x_0\right)=0& , & \X \in \mathbb{R}^n.
%		\end{array}\right.
%	\end{equation}
\end{proposition}

Proof:
	Starting from \eqref{DiracRieszFellerFourier}, we recall that the Cauchy problem (\ref{CauchySpaceFractionalDirac}) in the Fourier domain is given by
	\begin{equation}%\tag{$CP_\alpha^+$}
		\label{CauchySpaceFractionalFourier}\left\{\begin{array}{lll} 
			\displaystyle 	\left(~\partial_{x_0}+\vert \XI \vert^\alpha h_\theta(\XI) ~\right)\widehat{\u}(\XI,x_0)=0 &,~\XI \in \BR^n\setminus \{\bm 0\}&,~x_0\in \mathbb{R} \setminus \{0\}
			\\[1ex]
			\displaystyle	\widehat{\u}(\XI,0)=\widehat{\f}(\XI)&,~\XI \in \BR^n\setminus \{\bm 0\}, \\[1ex]
			\displaystyle	\lim_{\vert x_0 \vert \rightarrow +\infty}\widehat{\u}(\XI,x_0)=0&,~\XI \in \BR^n\setminus \{\bm 0\},
		\end{array}\right.
	\end{equation}
	where 
	\begin{center}
		$\displaystyle h_\theta(\XI)=\chi_-(\XI)+e^{-i{\theta\pi}}\chi_+(\XI)$, $\widehat{\f}(\XI)=(\mathcal{F}\f)(\XI)$ \& $\widehat{\u}(\XI,x_0)=(\mathcal{F}\u(\cdot, x_0))(\XI)$.
	\end{center}
	
Then, from the definition \eqref{uAnsatz} of $\u$ and from the pseudo-differential representation \eqref{RieszFellerHardy} of $\D_\theta^\alpha$, one easily finds that $\widehat{\u}(\XI,x_0)$ solves \eqref{CauchySpaceFractionalFourier} if, and only if, $\chi_-(\XI)\widehat{\u}(\XI,x_0)$ and $\chi_+(\XI)\widehat{\u}(\XI,x_0)$ solve

	\begin{eqnarray*}
	\label{CauchySpaceFractionalMinus} \left\{\begin{array}{lll} 
		\left(~\partial_{x_0}+\vert\XI\vert^\alpha~\right) \left(~\chi_-(\XI)\widehat{\u}(\XI,x_0)~\right)=0, &,~\XI \in \BR^n\setminus \{\bm 0\}&,~x_0>0
		\\ \ \\
		\displaystyle \lim_{x_0 \rightarrow 0^+}\chi_-(\XI)\widehat{\u}(\XI,x_0)=\chi_-(\XI)\widehat{\f}(\XI)&,~\XI \in \BR^n\setminus \{\bm 0\} \\ \ \\
		\displaystyle \lim_{x_0 \rightarrow +\infty}\chi_-(\XI)\widehat{\u}(\XI,x_0)=\chi_-(\XI)\widehat{\f}(\XI)& ~\XI \in \BR^n\setminus \{\bm 0\}
	\end{array}\right.
\end{eqnarray*}
and 	\begin{eqnarray*}
	\label{CauchySpaceFractionalPlus} \left\{\begin{array}{lll} 
		\left(~\partial_{x_0}+e^{-i\pi \theta}\vert\XI\vert^\alpha~\right) \left(~\chi_+(\XI)\widehat{\u}\left(\XI,x_0\right)~\right)=0, &,~\XI \in \BR^n\setminus \{\bm 0\}&,~x_0<0
		\\ \ \\
		\displaystyle \lim_{x_0 \rightarrow 0^-}\chi_+(\XI)\widehat{\u}(\XI,x_0)=\chi_+(\XI)\widehat{\f}(\XI)&,~\XI \in \BR^n\setminus \{\bm 0\} \\ \ \\
		\displaystyle \lim_{x_0 \rightarrow -\infty}\chi_+(\XI)\widehat{\u}(\XI,x_0)=\chi_+(\XI)\widehat{\f}(\XI)&,~\XI \in \BR^n\setminus \{\bm 0\}.
	\end{array}\right.
\end{eqnarray*}

%	and 	\begin{eqnarray}
%		\label{CauchySpaceFractionalPlus} \left\{\begin{array}{lll} 
%			\left(~\partial_{x_0}+e^{-i\pi \theta}\vert\XI\vert^\alpha~\right) \left(~\chi_+(\XI)\widehat{\u}\left(\XI,x_0\right)~\right)=0, &,~\XI \in \BR^n\setminus \{\bm 0\}&,~x_0<0
%			\\ \ \\
%			\displaystyle \lim_{x_0 \rightarrow 0^-}\chi_+(\XI)\widehat{\u}(\XI,x_0)=\chi_+(\XI)\widehat{\f}(\XI)&,~\XI \in \BR^n\setminus \{\bm 0\} \\ \ \\
%			\displaystyle \lim_{x_0 \rightarrow -\infty}\chi_+(\XI)\widehat{\u}(\XI,x_0)=\chi_+(\XI)\widehat{\f}(\XI)&,~\XI \in \BR^n\setminus \{\bm 0\}.
%		\end{array}\right.
%	\end{eqnarray}
	
	Using the method of characteristics, we find that $\chi_-(\XI)\widehat{\u}(\XI,x_0)$ and $\chi_+(\XI)\widehat{\u}(\XI,x_0)$ are uniquely determined by
\begin{eqnarray*}
	\chi_-(\XI)\widehat{\u}(\XI,x_0)=e^{-x_0\vert\XI\vert^\alpha}\widehat{\f_+}(\XI) &\mbox{\&}& \chi_+(\XI)\widehat{\u}(\XI,x_0)=e^{-x_0e^{-i\pi \theta}\vert\XI\vert^\alpha}\widehat{\f_-}(\XI),
\end{eqnarray*}
where $\f_\pm=\frac{1}{2}\left(\f\pm\mathcal{H}\f\right)$.

	Next, using the Fourier inversion formula, one finds that they coincide with \eqref{uPlus} (case of $\beta=0$) and \eqref{uMinus} (case of $\beta=\theta$), respectively. That is,
	\begin{center}
		$\mathcal{F}^{-1}\left(~\chi_-(\XI)\widehat{\u}(\XI,x_0)~\right)=\u_+(\X,x_0)$ \& $\mathcal{F}^{-1}\left(~\chi_+(\XI)\widehat{\u}(\XI,x_0)~\right)=\u_-(\X,x_0)$.
	\end{center} 
	
In addition, we emphasize that 
\begin{center}
	$\displaystyle \lim_{x_0\rightarrow 0^+}\chi_-(\XI)\widehat{\u}(\XI,x_0)=\chi_-(\XI)\widehat{\f}(\XI)$ \& $\displaystyle \lim_{x_0\rightarrow +\infty}\chi_-(\XI)\widehat{\u}(\XI,x_0)=0$
\end{center} for almost all $\XI\in\BR^n\in \setminus \{\bm 0\}$,
is derived from the exponential decay of $e^{-x_0\vert \XI\vert^\alpha}$, for values of $x_0>0$, 
while for values of $x_0<0$ and $\vert 1-\theta \vert < \frac{\alpha}{2}$, the limit conditions are
\begin{center}
	$\displaystyle \lim_{x_0\rightarrow 0^-}\chi_+(\XI)\widehat{\u}(\XI,x_0)=\chi_+(\XI)\widehat{\f}(\XI)$ \& $\displaystyle \lim_{x_0\rightarrow -\infty}\chi_+(\XI)\widehat{\u}(\XI,x_0)=0$
\end{center}
for almost all $\XI\in\BR^n$, resulting in
$\displaystyle \left\vert~e^{-x_0e^{-i\pi \theta}\left\vert\XI\right\vert^\alpha}~\right\vert = e^{-x_0\cos(\pi \theta)\vert\XI\vert^\alpha}. $

Thus, $\u_+$ and $\u_-$ are the unique solutions of (\ref{CauchySpaceFractionalHplus}).
%and (\ref{CauchySpaceFractionalHminus}), respectively. 
$\square$

\subsection{Cauchy type kernels $E_{\alpha,n}^\pm$}\label{subsec43}

The characterization of solutions to the Cauchy problem \eqref{CauchySpaceFractionalDirac}, as established in Proposition \ref{CauchyProblemHardy}, provides a groundbreaking framework for studying fractional analogues of the Cauchy kernel.\\[-3ex] 
\begin{eqnarray} \label{CauchyKernel} \E(\X,t)=\dfrac{1}{2\omega_n}\frac{t+\X}{\left(t^2+\vert \X \vert^2\right){\frac{n+1}{2}}}, &(\X,t)\in \BR^{n+1}_+ \end{eqnarray} 
building on Feller's influential technique \cite{Feller62} for stable distributions (see \cite[Chapter VI]{Feller71} for an overview). Here and elsewhere, the constant $\displaystyle \omega_n=\dfrac{\pi^{\frac{n+1}{2}}}{\Gamma\left(\frac{n+1}{2}\right)}$, which represents the measure of the unit sphere $\mathbb{S}^n$ in $\BR^{n+1}$.

Following the line of reasoning considered by the second author in \cite{Faustino21}, let us define it by \\[-3ex]
\begin{eqnarray}
	\label{CauchKernel-alpha}
%\hspace*{1cm}
	\E^{\pm}_{\alpha,n}(\X,z)=\frac{1}{2}\left(~{\bm K}_{\alpha,n}(\X,z)\pm \mathcal{H}{\bm K}_{\alpha,n}(\X,z)~\right),~\X\in\BR^n,~\Re(z)>0,
\end{eqnarray}
where
${\bm K}_{\alpha,n}(\X,z)$ denotes the {\textit radially symmetric} kernel function given by 
	\begin{equation}%\tag{$\mathcal{K}_{\alpha,n}$}
	\label{PoissonKalpha}{\bm K}_{\alpha,n}(\X,z)=\frac{1}{(2\pi)^{n}}\int_{\BR^n}e^{-z\vert \XI \vert^\alpha} e^{i\langle \X,\XI\rangle}\dxi ,~~\X\in \BR^n,~~\Re(z)>0.
\end{equation}

When $\alpha=1$, the set of the identities
\begin{eqnarray}
	\label{CauchyKernelPM}
	\begin{cases}
		\E^{+}_{1,n}(\X,x_0)=\E(\X,x_0) , & (\X,x_0)\in \BR^{n+1}_+ \\[1ex]
		\E^{-}_{1,n}(\X,x_0)=-\E(\X,-x_0) , & (\X,x_0)\in \BR^{n+1}_-
	\end{cases}
\end{eqnarray}
%\begin{eqnarray}
%	\label{PoissonK1}
%	{\bm P}(\X)=\dfrac{\Gamma\left(\frac{n+1}{2}\right)}{\pi^{\frac{n+1}{2}}}~\frac{1}{\left(1+\vert \X \vert^2\right)^{\frac{n+1}{2}}}, & \X \in \BR^n
%\end{eqnarray}
result from the Fourier inversion formulas (cf.~\cite[Theorem 1.14]{SteinWeiss72}) 
\begin{eqnarray*}
	\begin{cases}
		\displaystyle	\tfrac{1}{2}\mathcal{F}^{-1}\left(e^{-t\vert \cdot \vert} \right)(\X)= \frac{1}{2\omega_{n}}\frac{t}{\left(~t^2+\vert \X\vert^2~\right)^{\frac{n+1}{2}}} \\[3ex] \displaystyle
		\tfrac{1}{2}\mathcal{F}^{-1}\left(e^{-t\vert \cdot \vert} \frac{-i\left(\cdot\right)}{\vert\cdot \vert} \right)(\X)=\frac{1}{2\omega_{n}}\frac{\X}{\left(~t^2+\vert \X \vert^2~\right)^{\frac{n+1}{2}}}
	\end{cases}, &(\X,t)\in  \BR^{n+1}_+.
\end{eqnarray*}

The Poisson kernel $\displaystyle {\bm K}_{1,n}(\X,\vert x_0\vert)=\E(\X,\vert x_0\vert)+\E(-\X,\vert x_0\vert)$, as encoded in the set of identities \eqref{CauchyKernelPM}, has also been considered in applications involving Brownian motion to describe a {\textit Cauchy probability distribution law} in the upper half space $\BR^{n+1}_+$, from the north pole of the $n-$dimensional unit sphere $\mathbb{S}^n$, as detailed in the papers by Letac \cite{Letac1986} and Kato-McCullagh \cite{KatoMcCullagh20} on conformal probability distributions. 
Since our current understanding of probability is not sufficient to discuss this on a technical level, we will rely solely on Feller's account \cite{Feller62} to establish a concrete connection between the Cauchy kernels \eqref{CauchKernel-alpha} and the one-parameter semigroup $\displaystyle \left\{e^{-x_0\D_\theta^\alpha}\right\}_{x_0\in \BR}$ generated by $-\D_\theta^\alpha$.

In a broad sense, for $z=te^{\frac{i\pi \gamma}{2}}$ we can interpret the Fourier transform of the kernel function ${\bm K}_{\alpha,n}$, defined by \eqref{PoissonKalpha}, as a stable distribution on $\BR^n$ with index $\alpha$ and skewness $\gamma$, where $0<\alpha\leq 1$ and $\vert \gamma\vert<\alpha$ (cf. ~\cite{BG1960,Mainardi01}). Namely, if  
\begin{eqnarray}
	\label{Characteristic}\widehat{\bm K}_{\alpha,n}(\XI,z)=e^{-z\vert \XI\vert^\alpha}, &\X\in \BR^n,~~\Re(z)>0,
\end{eqnarray}
represents the characteristic function of a {\textit symmetric stable process} $\left\{{\bm X}(z)~:~\Re(z)> 0\right\}$, it turns out that $\widehat{\bm K}_{\alpha,n}\left(\XI,te^{\frac{i\pi\gamma}{2}}\right)$ solves the Cauchy problem
\begin{eqnarray}
	\label{FundamentalSolutionFourier}
	\begin{cases}
		\partial_{t} \widehat{\bm K}_{\alpha,n}\left(\XI,te^{\frac{i\pi\gamma}{2}}\right)+e^{\frac{i\pi\gamma}{2}}\vert \XI\vert^\alpha \widehat{\bm K}_{\alpha,n}\left(\XI,te^{\frac{i\pi\gamma}{2}}\right)=0 &,~(\XI,t)\in \BR^{n+1}_+ \\[2ex]
		\displaystyle \lim_{t\rightarrow 0^+}	\widehat{\bm K}_{\alpha,n}\left(\XI,te^{\frac{i\pi\gamma}{2}}\right)=1 &,~\XI\in \BR^n.
	\end{cases}
\end{eqnarray}

Thus, by taking the Fourier inverse on both sides of \eqref{FundamentalSolutionFourier}, it turns out that \eqref{PoissonKalpha} provides an integral representation of the fundamental solution of the space-fractional operator $\partial_t+e^{\frac{i\pi\gamma}{2}}(-\Delta)^{\frac{\alpha}{2}}$. That is, ${\bm K}_{\alpha,n}\left(\X,te^{\frac{i\pi\gamma}{2}}\right)$ solves
\begin{eqnarray}
	\label{FundamentalSolution}
	\begin{cases}
		\partial_{t} {\bm K}_{\alpha,n}\left(\X,te^{\frac{i\pi\gamma}{2}}\right)+e^{\frac{i\pi\gamma}{2}}(-\Delta)^{\frac{\alpha}{2}} {\bm K}_{\alpha,n}\left(\X,te^{\frac{i\pi\gamma}{2}}\right)=0 &,~(\X,t)\in \BR^{n+1}_+ \\[2ex]
		\displaystyle \lim_{t\rightarrow 0^+}{\bm K}_{\alpha,n}\left(\X,te^{\frac{i\pi\gamma}{2}}\right)=\delta(\X) &,~\X\in \BR^n,
	\end{cases}
\end{eqnarray}
where $\delta$ is the delta function on $\BR^n$.

The previous identity allows us to represent, in the distributional sense, the analytic semigroup of L\'evy-Feller type, $\left\{\exp\left(-te^{\frac{i\pi \gamma}{2}}(-\Delta)^{\frac{\alpha}{2}}\right)\right\}_{t\geq 0}$, with the parameters $0<\alpha\leq 1$ and $\vert \gamma\vert<\alpha$,
by the convolution formula
\begin{eqnarray}
	\label{LevyFellerSemigroup}
	\exp\left(-te^{\frac{i\pi \gamma}{2}}(-\Delta)^{\frac{\alpha}{2}}\right)\f={\mathbf K}_{\alpha,n}\left(\cdot,te^{\frac{i\pi \gamma}{2}}\right)*\f.
\end{eqnarray}

Based on this formula, we can obtain a similar characterization for $\displaystyle \left\{e^{-x_0\D_\theta^\alpha}\right\}_{x_0\in \BR}$. This leads to the following proposition:

\begin{proposition}\label{CauchyProblemHardy-Convolution}
Let $\E^\pm_{\alpha,n}$ denote the Cauchy kernels defined by \eqref{CauchKernel-alpha}.	Then, for each $\f\in L^p(\BR^n;\cl_{0,n})$, the Clifford valued function 
\begin{eqnarray}
	\label{FourierInversionU}
	\u(\X,x_0)= \begin{cases}
		\left(\E_{\alpha,n}^+\left(\cdot,x_0\right)*\f\right)(\X)&,~(\X,x_0)\in \BR^{n+1}_+ \\[1ex] \left(\E_{\alpha,n}^-\left(\cdot,x_0e^{-i\pi\theta}\right)*\f\right)(\X)&,~(\X,x_0)\in \BR^{n+1}_- \\[1ex]
		\f(\X)&,~(\X,x_0)\in \BR^{n}\times \{0\}
	\end{cases}
\end{eqnarray}
solves the Cauchy problem (\ref{CauchySpaceFractionalDirac}) if and only if
\begin{eqnarray}
	\label{Fpm} \u(\X,x_0)=e^{-x_0\D_\theta^\alpha}{\bm F}(\X,x_0)&,~(\X,x_0)\in \BR^{n+1}&,
\end{eqnarray}
where
\begin{eqnarray}\label{Fpm-Hardy}
	{\bm F}(\X,x_0)=\begin{cases}
		\displaystyle \tfrac{1}{2}(\f(\X)+\mathcal{H}\f(\X))&,~(\X,x_0)\in \BR^{n+1}_+ \\[1ex] \displaystyle \tfrac{1}{2}(\f(\X)-\mathcal{H}\f(\X))&,~(\X,x_0)\in \BR^{n+1}_- \\[1ex]
		\f(\X) &,~(\X,x_0)\in \BR^{n}\times \{0\}.
	\end{cases}
\end{eqnarray}
\end{proposition}

Proof:
From the convolution formula
\begin{eqnarray*}
	\left({\bm K}_{\alpha,n}\left(\cdot,z\right)*\f_\pm\right)(\X)=\frac{1}{(2\pi)^n}\int_{\BR^n}e^{-z\vert \XI\vert^{\alpha}}\chi_{\mp}(\XI)\widehat{\f}(\XI) e^{i\langle \X,\XI\rangle}\dxi, &\Re(z)>0,
\end{eqnarray*}
with $\displaystyle \f_\pm=\frac{1}{2}\left(\f\pm \mathcal{H}\f\right)$, it is clear that the function $\u$, defined by \eqref{FourierInversionU}, can be rewritten in terms of 
\begin{center}
	$\u_+\in \displaystyle H^p_{\alpha,0}\left(\BR^{n+1}_+;\cl_{0,n}\right)$ and $\u_-\in H^p_{\alpha,\theta}\left(\BR^{n+1}_-;\cl_{0,n}\right)$.
\end{center}

Namely, the proof that \eqref{FourierInversionU} solves the Cauchy problem \eqref{CauchySpaceFractionalDirac} follows from the fact that the functions $\u_\pm$, uniquely determined by
	\begin{eqnarray*}
		\begin{array}{llll}
			\u_+(\X,x_0)&=&	\left(\E_{\alpha,n}^+\left(\cdot,x_0\right)*\f\right)(\X)&,~(\X,x_0)\in \BR^{n+1}_+ \\[1ex] \u_-(\X,x_0)&=&\left(\E_{\alpha,n}^-\left(\cdot,x_0e^{-i\pi\theta}\right)*\f\right)(\X)&,~(\X,x_0)\in \BR^{n+1}_-,
		\end{array}
	\end{eqnarray*}
are under the conditions of Proposition \ref{CauchyProblemHardy}.
Thus, to show that $\u$ satisfies \eqref{Fpm}, it suffices to prove that
\begin{eqnarray}
	\label{EPlusMinus} \begin{array}{lll}
		\displaystyle e^{-x_0\vert \XI\vert^\alpha h_\theta(\XI)}\chi_-(\XI)=e^{-x_0\vert \XI\vert^\alpha}\chi_-(\XI), & x_0>0 \\[1ex]
		e^{-x_0\vert \XI\vert^\alpha h_\theta(\XI)}\chi_+(\XI) =e^{-x_0e^{-i\pi \theta}\vert \XI\vert^\alpha}\chi_+(\XI) , & x_0<0,
	\end{array}
\end{eqnarray}
where $\displaystyle h_\theta(\XI)=\chi_-(\XI)+e^{-i{\theta\pi}}\chi_+(\XI)$ $(~\XI \in \BR^n\setminus \{\bm 0\}~)$.
	
To proceed, recall that
$$
e^{-x_0\vert \XI\vert^\alpha h_\theta(\XI)}=\sum_{k=0}^\infty \frac{(-x_0)^k}{k!}\vert \XI\vert^{\alpha k} \left(~\chi_-(\XI)+e^{-i{\theta\pi}}\chi_+(\XI)~\right)^k
$$ 
represents the power series expansion of $\displaystyle e^{-x_0\vert \XI\vert^\alpha h_\theta(\XI)}$. Then, using the proof strategy outlined in \cite[Theorem 1]{Faustino21}, one can easily deduce the set of properties \eqref{ProjectionOpProperties} associated with the idempotents $\chi_\pm(\XI)$
	\begin{eqnarray*}
		\left(~\chi_-(\XI)+e^{-i\pi \theta}\chi_+(\XI)~\right)^k=\chi_-(\XI)+e^{-i\pi \theta k}\chi_+(\XI), &\forall~k\in \mathbb{N}.
	\end{eqnarray*}
	
	Thus, $e^{-x_0\vert \XI\vert^\alpha h_\theta(\XI)}$ simplifies to 
	\begin{eqnarray*}
		e^{-x_0\vert\XI\vert^\alpha h_\theta(\XI)}
		&=& \sum_{k=0}^\infty \frac{(-x_0\vert \XI\vert^{\alpha})^k}{k!}\chi_-(\XI)+\sum_{k=0}^\infty \frac{(-x_0e^{-i{\theta\pi}}\vert \XI\vert^{\alpha})^k}{k!}\chi_+(\XI)\\
		&=&e^{-x_0\vert \XI\vert^\alpha}\chi_-(\XI)+e^{-x_0e^{-i\pi \theta}\vert \XI\vert^\alpha}\chi_+(\XI),
	\end{eqnarray*}
	proving that $e^{-x_0\vert\XI\vert^\alpha h_\theta(\XI)}\chi_-(\XI)$ and $e^{-x_0\vert\XI\vert^\alpha h_\theta(\XI)}\chi_+(\XI)$ are given by                                                                                                                                                                                                                                                                                                                                                                                                                                                                                                                                                                                                                                                                                                                                                                                                                                                                                                                                                                                                                                                                                                                                                                                                                                                                                                                                                                                                                                                                                                                                                                                                                                                                                                                                                                          \eqref{EPlusMinus}. This completes our proof.
$\square$

%%%% @nelson - redundant remark
%\begin{remark}
%	From the proof of Proposition \ref{CauchyProblemHardy-Convolution}, it turns out that
%	\begin{equation*}
	%		\label{HypercomplexExp}e^{-x_0\vert \XI\vert^\alpha h_\theta(\XI)}e^{i\langle \X,\XI \rangle}=\left(e^{-x_0\vert\XI\vert^\alpha}\chi_-(\XI)+e^{-x_0e^{-i\pi \theta}\vert\XI\vert^\alpha}\chi_+(\XI)\right)e^{i\langle \X,\XI \rangle}
	%	\end{equation*}
%	may be interpreted as the hypercomplex extension of the Fourier kernel $e^{i\langle \X,\XI \rangle}$.
%\end{remark}

\section{Paley-Wiener type theorems}\label{sec6}

%%% @nelson
%In this section, we prove a series of Paley-Wiener equivalences that establish a tangible correspondence between the growth behavior of the sequences $\left(~\left(~\D_\theta^\alpha~\right)^k \f_\pm~\right)_{k \in \BN_0}$, the growth behavior of the solutions of the Cauchy problem \eqref{CauchySpaceFractionalDirac}, and the support of $\widehat{\f}=\mathcal{F}\f$, in the case where  $\operatorname{supp}\widehat{\f}\subseteq \overline{B(0,R)}$.

%In Subsection \ref{subsec61} we start with the proof of Theorem \ref{RealPaleyWiener-BumpFunction} for Schwartz bump functions. Such a proof strategy mimics the approach considered by Andersen and Jeu in \cite{AndersenJeu10}. Afterwards, we will use Theorem \ref{RealPaleyWiener-BumpFunction} to prove Theorem \ref{PaleyWienerEquivalences} in Subsection \ref{subsec62}. This result corresponds to a real Paley-Wiener type theorem in the context of hypercomplex variables.

%To get in tocuOur generalization of Bernstein type spaces takes into account the solutions of the Cauchy problem \eqref{CauchySpaceFractionalDirac} of exponential type $R^\alpha$, which is closely related with the approach of Pesenson and Zayed on abstract Paley-Wiener spaces (cf.~\cite{Pesenson01,PesensonZayed09}). That gives in turn a faithful generalization for the Bernstein spaces, appearing in the proof of \cite[Theorem 5.1]{FranklinHoganLarkin17} and elsewhere.

\subsection{The radially symmetric case}\label{subsec61}

First, for each $k,\ell \in \BN_0$ and $\X \in \BR^n$, let us take a closer look at the integral representation of a \textit{radially symmetric function} $\PSI\in \mathcal{S}(\BR^n;\cl_{0,n})$:	\begin{eqnarray}\label{FourierInversionSchwarz}
	\vert \X \vert^{2\ell}\left(~(-\Delta)^{\frac{\alpha}{2}}~\right)^k\PSI(\X)~ =\frac{1 }{(2\pi)^n}\int_{\BR^n} (-\Delta_{\XI})^\ell\left(\vert \XI \vert^{\alpha k} \widehat{\PSI}(\XI)\right)~e^{i\langle \X,\XI\rangle}\dxi,
\end{eqnarray}
where $\Delta_{\XI}$ is the Laplace operator on the variable $\XI$. In terms of the polar coordinates $\rho=\vert \XI\vert$ and $\displaystyle {\bm \omega}=\frac{\XI}{\vert \XI \vert}$ ($\XI\in \BR^n\setminus \{\bm 0\}$), we get the split
\begin{eqnarray}
	\Delta_{\XI}=\frac{d^2}{d\rho^2}+\frac{n-1}{\rho}\frac{d}{d\rho}+\frac{1}{\rho^2}\Delta_{\mathbb{S}^{n-1}},
\end{eqnarray}
where $\Delta_{\mathbb{S}^{n-1}}$ stands for the Laplace-Beltrami operator on the $(n-1)$-dimensional sphere $\mathbb{S}^{n-1}$.

An important property of $\Delta_{\XI}$ is that it maps {\textit radially symmetric functions} to {\textit radially symmetric functions}. 
Furthermore, based on the following Fourier inversion formula for {\textit radially symmetric functions} $\varphi$: 
\begin{eqnarray}
	\label{FourierInversionRadial}
	\int_{\BR^n} \varphi(\vert\XI\vert)e^{i \langle \X, \XI \rangle} d\XI =\frac{(2\pi)^{\frac{n}{2}}}{\vert\X\vert^{\frac{n-2}{2}}} \int_{0}^{\infty} \varphi(\rho) \rho^{\frac{n}{2}}J_{\frac{n}{2}-1}(\rho\vert\X\vert)d\rho
\end{eqnarray}
there holds the following equivalent representation for \eqref{FourierInversionSchwarz}:  
\begin{eqnarray}
	\label{FourierInversionSchwarz-Hankel}
	\begin{array}{lll}
		\vert \X \vert^{2\ell}\left(~(-\Delta)^{\frac{\alpha}{2}}~\right)^k\PSI(\X)=\\=\displaystyle \frac{1}{(2\pi)^{\frac{n}{2}}}\int_{0}^\infty \left(-\widetilde{\Delta}_{\frac{n}{2}-1}\right)^\ell\left(\rho^{\alpha k} \widehat{\PSI}(\rho)\right)~\frac{J_{\frac{n}{2}-1}(\rho\vert \X \vert)}{(\rho\vert \X \vert)^{\frac{n}{2}-1}}~\rho^{n-1}d\rho,
	\end{array}
\end{eqnarray}
$J_\lambda$ stands for the Bessel function of order $\lambda$ and
\begin{eqnarray}
	\label{BesselOp} \displaystyle \widetilde{\Delta}_{\lambda}:=\frac{d^2}{d\rho^2}+\frac{2\lambda+1}{\rho}\frac{d}{d\rho}
\end{eqnarray} for the Bessel operator of order $\lambda$, with $\lambda>-\frac{1}{2}$.

%%% @nelson 
%We recall that the normalized Bessel functions $\widetilde{J}_{\nu}(z)=z^{-\nu}{{J}_{\nu}(z)}$, which appear on the right side of \eqref{FourierInversionSchwarz-Hankel},
%are eigenfunctions of the Bessel operator \eqref{BesselOp}. Furthermore, from \eqref{FourierInversionRadial} and
%$$
%\int_{0}^{1} J_{\mu}(z\rho)\rho^{\mu+1}(1-\rho^2)^\nu d\rho=\frac{\Gamma(\nu+1)2^\nu}{z^{\nu+1}}J_{\mu+\nu+1}(z),
%$$
%which holds for $\Re(\mu)>-\frac{1}{2}$, $\Re(\nu)>-1$ and $z>0$ (see \textcolor{red}{\ldots}),
%we conclude that the Fourier transform of $\widetilde{J}_{\lambda+\frac{n}{2}-1}\left(R\vert \cdot\vert\right)$ is a real-valued bump function with support $\overline{B(0,R)}$, since
%\begin{eqnarray}
%	\label{BumpFunction-Bessel}
%	~\frac{J_{\lambda+\frac{n}{2}-1}(R\vert \X \vert)}{\left(~R\vert \X \vert~\right)^{\lambda+ \frac{n}{2}-1}}=\frac{1}{R^{\frac{n}{2}+1}}(\mathcal{F}^{-1}{\bm m}_\lambda)(\X), & \X\in \BR^n,
%\end{eqnarray}
%where $\left({\bm m}_\lambda\right)_{\lambda>0}$ denotes the family of functions defined by
%\begin{eqnarray}
%	\label{BumpFunction}
%	{\bm m}_\lambda(\XI)=\begin{cases}\displaystyle \frac{1}{\Gamma(\lambda)2^{\lambda-1}}
%		\left(1-\frac{\vert \XI\vert^2}{R^2}\right)^{\lambda} &,~\vert\XI \vert \leq R \\ \ \\
%		0&,~\vert \XI\vert > R.
%	\end{cases}
%\end{eqnarray}

Thus, the integral representation \eqref{FourierInversionSchwarz-Hankel} provides an effective approach to deriving a {\textit real Paley-Wiener theorem} by real variable methods. The following theorem can be interpreted as a reformulation of \cite[Theorem 3]{Andersen06} for Clifford valued functions that are {\textit radially symmetric}.

\begin{theorem}\label{RealPaleyWiener-BumpFunction}
	Let $\PSI\in\mathcal{S}(\BR^n;\cl_{0,n})$ be a {\textit radially symmetric function}. Then the following statements are equivalent:
\begin{enumerate}
	\item[\textbf (a)] $\operatorname{supp}\widehat{\PSI}\subseteq \overline{B(0,R)}$.
	\item[\textbf (b)] For each $m\in \mathbb{N}_0$, there exists a constant $0<\Lambda_{m,n}<\infty$ such that for all $k\in \BN_0$, satisfying the condition $\alpha k>2m+1-n$, we have
	\begin{eqnarray}%\tag{$\mathcal{B}_{\alpha k}$}
		\label{PW-BumpIneq}
		% \begin{array}{lll}
			\displaystyle \sup_{\X \in \BR^n}~\left\vert~\left(1+\vert \X \vert^2\right)^m\left(~(-\Delta)^{\frac{\alpha}{2}}~\right)^k\PSI(\X)~\right\vert_0\leq \Lambda_{m,n}~R^{\alpha k}.
			% \end{array}
	\end{eqnarray}
\end{enumerate}
\end{theorem}

\textbf{Proof of  (a) $\Longrightarrow$ (b)}

First, let us take a closer look at the sequence of functions $\left(\widetilde{\Phi}_{\alpha k,\ell}\right)_{\ell\in \BN_0}$, defined by
	\begin{eqnarray*}
		\widetilde{\Phi}_{\alpha k,\ell}(\rho)=\left(-\widetilde{\Delta}_{\frac{n}{2}-1}\right)^\ell\left(\rho^{\alpha k} \widehat{\PSI}(\rho)\right)&,~\rho>0.
	\end{eqnarray*}
	
For $\ell=0$ there is nothing to prove, while for $\ell=1$ the closed formula
	$$
	\widetilde{\Phi}_{\alpha k,1}(\rho)=- s_1(\alpha k,n) \rho^{\alpha k-2}\left(\widehat{\PSI}(\rho)+\frac{1}{s_1(\alpha k,n)}\rho^2\widetilde{\Delta}_{\alpha k+\frac{n}{2}-1}\widehat{\PSI}(\rho)\right),
	$$
	with $s_1(\alpha k,n)=\alpha k(\alpha k+(n-2))$,
	yields from \eqref{BesselOp} and from the set of relations
	\begin{eqnarray*}
		\frac{1}{\rho}\frac{d}{d\rho}\left(\rho^{\alpha k} \widehat{\PSI}(\rho) \right)&=&\alpha k\rho^{\alpha k-2}\widehat{\PSI}(\rho)+\rho^{\alpha k-1}\frac{d}{d\rho}\widehat{\PSI}(\rho), \\
		\frac{d^2}{d\rho^2}\left(\rho^{\alpha k} \widehat{\PSI}(\rho) \right)&=&\alpha k(\alpha k-1)\rho^{\alpha k-2}\widehat{\PSI}(\rho)
		+ 2\alpha k\rho^{\alpha k-1}\frac{d}{d\rho}\widehat{\PSI}(\rho)+\rho^{\alpha k}\frac{d^2}{d\rho^2}\widehat{\PSI}(\rho).
	\end{eqnarray*}
	
Consequently, the assumption that $\operatorname{supp}\widehat{\PSI} \subseteq \overline{B(0,R)}$, combined with the fact that $\widehat{\PSI}$ can be represented in terms of \eqref{FourierInversionRadial}, establishes the existence of two constants $\lambda_0 > 0$ and $\lambda_1 > 0$ such that
	\begin{eqnarray*}
		\vert~\widehat{\PSI}(\XI)~\vert_0\leq \lambda_0 &\mbox{\&}& \left\vert~\widetilde{\Delta}_{\alpha k+\frac{n}{2}-1}\widehat{\PSI}(\XI)~\right\vert_0 \leq \lambda_1
	\end{eqnarray*}
	hold for every $\vert \XI\vert<R$. Hence,
	\begin{eqnarray*}
		\vert \widetilde{\Phi}_{\alpha k,1}(\rho)\vert_0\leq \lambda(\alpha k+(n-1))^2 \rho^{\alpha k-2} &,~ 0<\rho\leq R,
	\end{eqnarray*}
	for some $\lambda>0$.
	
	Then, by performing an induction argument over $\ell\in \BN_0$, we are able to construct a sequence of {\textit radially symmetric functions} $\left({\PSI_\ell}\right)_{\ell \in \BN_0}$ such that
	\begin{eqnarray}
		\label{EstimateInduction}
		\left\vert~\widetilde{\Phi}_{\alpha k,\ell}(\rho)~\right\vert_0\leq (\alpha k+(n-1))^{2\ell} \rho^{\alpha k-2\ell}\left\vert~ \widehat{\PSI_\ell}(\rho)~\right\vert_0, & 0<\rho\leq R,
	\end{eqnarray}
	is true for any $\ell\in \BN_0$, where $\operatorname{supp}\widehat{\PSI_\ell}\subseteq \overline{B(0,R)}$.
	
	Further, 
	from the estimate \eqref{EstimateInduction} and from \eqref{FourierInversionSchwarz-Hankel} we get
	\begin{eqnarray*}
		\label{FourierInversionSchwarz-Ineq}
		%	\begin{array}{lll}
			\vert \X \vert^{2\ell}\left\vert~ \left(~(-\Delta)^{\frac{\alpha}{2}}~\right)^k\PSI(\X)\right\vert_0 \leq 
			\frac{1}{(2\pi)^{\frac{n}{2}}}\int_{0}^R	\left\vert~\widetilde{\Phi}_{\alpha k,\ell}(\rho)~\right\vert_0 ~\rho^{n-1}d\rho,
			\\
			\leq \displaystyle \frac{\Gamma\left(\frac{n}{2}\right)}{\pi^{\frac{n}{2}}}\|\widehat{\PSI_\ell}\|_\infty~(\alpha k+(n-1))^{2\ell} \int_{0}^R ~\rho^{\alpha k-2\ell+n-1}d\rho \\
			\\ \leq \frac{\Gamma\left(\frac{n}{2}\right)}{\pi^{\frac{n}{2}}}\|\widehat{\PSI_\ell}\|_\infty (\alpha k+(n-1))^{2\ell} R^{\alpha k-2\ell+n}.
			%	\end{array}
	\end{eqnarray*}
	
Therefore, by a simple calculation based on the binomial identity
	\begin{eqnarray*}
		%	\label{BinomialExpansionLaplacian}
		\left(~1+\vert \X\vert^2~\right)^m\left(~(-\Delta)^\alpha~\right)^k\PSI(\X)=\sum_{\ell=0}^m\left(\begin{array}{c} m \\
			\ell 
		\end{array}\right)\vert \X\vert^{2\ell}\left(~(-\Delta)^\alpha~\right)^k\PSI(\X),
	\end{eqnarray*}	
	we conclude that
$$	\left(~1+\vert \X \vert^{2}~\right)^m\left\vert~ \left(~(-\Delta)^{\frac{\alpha}{2}}~\right)^k\PSI(\X)\right\vert_0 \leq \Lambda_{m,n} R^{\alpha k+n}\left(1+\frac{(\alpha k+(n-1))^2}{R^2}\right)^m,$$
	  where
	$$\Lambda_{m,n}=\frac{\Gamma\left(\frac{n}{2}\right)}{\pi^{\frac{n}{2}}}~\max_{0\leq \ell\leq m}\|\widehat{\PSI_\ell}\|_\infty.$$ 

Furthermore, by noticing that the sequence $\left(\phi_k(R)\right)_{k\in \BN_0}$, given by
\begin{eqnarray*}
	\phi_k(R)=R^{\alpha k+n}\left(1+\frac{(\alpha k+(n-1))^2}{R^2}\right)^m,~~k\in \BN_0
\end{eqnarray*}
satisfies the $\limsup-$condition $\displaystyle \limsup_{k\rightarrow \infty}\left(\phi_k(R)\right)^{\frac{1}{k}}=R^\alpha$, we can deduce from the properties of $\limsup$ that the inequality \eqref{PW-BumpIneq} holds for the constant $\Lambda_{m,n}$, as given above. $\square$
\\[2ex]
\textbf{Proof of  (b) $\Longrightarrow$ (a)} Suppose $\vert \XI\vert>R+\varepsilon$ holds for any $\varepsilon>0$. Choose $m \in \BN_0$ so that $\X \mapsto \left(~1+\vert\X\vert^2~\right)^{-m}$ defines a $L^1-$function.
This allows us to rewrite $\left\|~\left(~(-\Delta)^{\frac{\alpha}{2}}~\right)^k\PSI~\right\|_1$ as follows (see Definition \ref{LpSpaces}):
	\begin{eqnarray*}
		\left\|~\left(~(-\Delta)^{\frac{\alpha}{2}}~\right)^k\PSI~\right\|_1
		=\int_{\BR^n} \left(~1+\vert \X\vert^2~\right)^{-m}~\left\vert~ \left(~1+\vert\X\vert^2~\right)^m\left(~(-\Delta)^{\frac{\alpha}{2}}~\right)^k\PSI~\right\vert_0 \dx,
	\end{eqnarray*}
	where $\vert \cdot \vert_0$ denotes the underlying norm of $\cl_{0,n}$.
	Furthermore, using the inequality~\eqref{PW-BumpIneq}, we get
	\begin{eqnarray*}
		\left\|~\left(~(-\Delta)^{\frac{\alpha}{2}}~\right)^k{\bm \psi}~\right\|_1 &\leq&  \left\|~\left(~1+\vert \cdot\vert^2~\right)^{-m}~\right\|_1
		 \sup_{\X \in \BR^n}\left\vert~ \left(~1+\vert \X \vert^2~\right)^m\left(~(-\Delta)^{\frac{\alpha}{2}}~\right)^k\PSI(\X)~\right\vert_0.
	\end{eqnarray*}
	
	Consequently, given that the statement \textbf{(b)} is true, the following set of inequalities holds for each $k \in \BN_0$:
	\begin{eqnarray*}
		(R+\varepsilon)^{\alpha k} \vert \widehat{\PSI}(\XI)\vert_0 \leq \vert \XI\vert^{\alpha k} \vert \widehat{\PSI}(\XI)\vert_0  
		\leq \left\|~\left(~(-\Delta)^{\frac{\alpha}{2}}~\right)^k \PSI~\right\|_1 
		\leq \Lambda_{m,n}~R^{\alpha k+n}\left\|~\left(~1+\vert \cdot \vert^2~\right)^{-m}~\right\|_1.
	\end{eqnarray*}
	
Furthermore, by observing that the sequence $\left(\phi_k(R,\varepsilon)\right)_{k\in \BN_0}$, given by
	\begin{eqnarray*}
\phi_k(R,\varepsilon)=\Lambda_{m,n}R^n\left(\frac{R}{R+\varepsilon}\right)^{\alpha k}&,~k\in \BN_0
	\end{eqnarray*}
	satisfies the limit condition $$\displaystyle \lim_{k\rightarrow \infty}\frac{\phi_{k+1}(R,\varepsilon)}{\phi_k(R,\varepsilon)}= \left(\frac{R}{R+\varepsilon}\right)^{\alpha},$$ we infer 
	\begin{eqnarray*}
		\left\vert~\widehat{\PSI}(\XI)~\right\vert_0\leq  \lim_{k\rightarrow \infty}\vartheta_{m,n}R^n\left(\frac{R}{R+\varepsilon}\right)^{\alpha k}\left\|~\left(~1+\vert \cdot \vert^2~\right)^{-m}~\right\|_1=0.
	\end{eqnarray*}
	So, $\widehat{\PSI}(\XI)=0$ for $\XI\not \in\operatorname{supp}\widehat{\PSI}$, as desired.
	
$\square$

\subsection{Real Paley-Wiener type theorems}\label{subsec62}

The following theorem, based on the philosophy of {\textit real Paley-Wiener type theorems}, establishes a concrete correspondence between the support of $\widehat{\f}=\mathcal{F}\f$, in the case where  $\operatorname{supp}\widehat{\f}\subseteq \overline{B(0,R)}$, the growth behavior of the sequences $\left(~\left(~\D_\theta^\alpha~\right)^k \f_\pm~\right)_{k \in \BN_0}$ and the growth behavior of the solutions of the Cauchy problem \eqref{CauchySpaceFractionalDirac}.
To prove this theorem using Theorem \ref{RealPaleyWiener-BumpFunction}, the following lemmas are required.

\begin{lemma}\label{Uryshon-Lemma}
	If $\operatorname{supp}\widehat{\f}\subseteq \overline{B(0,R)}$, then for any $\varepsilon>0$ there exists a real-valued function $\PSI$ such that
	\begin{description}\itemsep0.5ex
		\item[\textbf (a)] $\widehat{\PSI}\in C^\infty_0\left(B(0,R+\varepsilon)\right)$;
		\item[\textbf (b)] $0\leq \widehat{\PSI}(\XI)\leq 1$, for all $\XI \in B(0,R+\varepsilon)$;
		\item[\textbf (c)] $\widehat{\PSI}(\XI)=1$, for all $\XI\in \overline{B(0,R)}$.
	\end{description}
\end{lemma}

Proof:
Since $\overline{B(0,R)}$ is a compact set, $B(0,R+\varepsilon)$ is an open set, and $\overline{B(0,R)}\subseteq B(0,R+\varepsilon)$, the proof of Lemma \ref{Uryshon-Lemma} follows directly from Uryshon's Lemma. For more details, see~\cite[Exercise 1.15]{LiebLoss2001}.
$\square$

\begin{lemma}\label{DiracHardyRn-Lemma}
	Let $L^{p,\pm}(\BR^n;\cl_{0,n})$ denote the function spaces \eqref{HardySpacesRn}. Then, for every $\f_+\in L^{p,+}(\BR^n;\cl_{0,n})$ and $\f_-\in L^{p,-}(\BR^n;\cl_{0,n})$, we find that
\begin{eqnarray}
	\label{Dpm}
	\begin{array}{lll}
		\left(~\D_\theta^\alpha~\right)^k \f_+(\X)&=&\left(~(-\Delta)^{\frac{\alpha}{2}}~\right)^k\f_+(\X) \\[1ex] \left(~\D_\theta^\alpha~\right)^k\f_- (\X)&=&e^{-i\pi k \theta}\left(~(-\Delta)^{\frac{\alpha}{2}}~\right)^k\f_-(\X)
	\end{array}&,~\X \in \BR^n,
\end{eqnarray}
for each $k\in \BN_0$.
\end{lemma}

Proof:
	Recall that the set of relations
	\begin{eqnarray*}
		\begin{array}{lll}
			\displaystyle \D_\theta^\alpha \f_+(\X)&=&(-\Delta)^{\frac{\alpha}{2}}\f_+(\X)  \\[1ex] \displaystyle \D_\theta^\alpha \f_-(\X)&=&e^{-i\pi \theta}(-\Delta)^{\frac{\alpha}{2}}\f_-(\X)
		\end{array}&,~\X \in \BR^n
	\end{eqnarray*}
follows from
	\eqref{RieszFellerHardy} and from the fact that $\frac{1}{2}(I\pm \mathcal{H})$ are idempotent operators. 
	
	Furthermore, by using inductive arguments, the proof of \eqref{Dpm} is then immediate.
$\square$

\begin{theorem}\label{PaleyWienerEquivalences}
	Let $\u$ be the solution of the Cauchy problem \eqref{CauchySpaceFractionalDirac} such that 
\begin{eqnarray*}
	\u(\cdot,0)=\f&\mbox{and}&\f\in L^p(\BR^n;\cl_{0,n}). 
\end{eqnarray*}
Then, for any $1< p<\infty$ and $R>0$, the following statements are equivalent whenever the closed formula of $\u$ is given by Proposition \ref{CauchyProblemHardy-Convolution}:
	\begin{description}
		\item[\textbf (A)]  $\operatorname{supp}\widehat{\f}\subseteq \overline{B\left(0,R\right)}$.
		\item[\textbf (B)] For each $k\in \BN_0$, there exists $C_+>0$ resp. $C_->0$ such that the so-called Bernstein inequalities
		\begin{eqnarray*}
			~\left\|~\left(~\D_\theta^\alpha~\right)^k\f_\pm~\right\|_p\leq C_\pm R^{\alpha k}\|~\f_\pm~\|_p.
		\end{eqnarray*}
		hold for $\displaystyle \f_+=\frac{1}{2}\left(\f+ \mathcal{H}\f\right)$ and $\displaystyle \f_-=\frac{1}{2}\left(\f-\mathcal{H}\f\right)$, respectively.
		\item[\textbf (C)] There exists $C>0$ such that 
		$$ ~\left\|~\u(\cdot,x_0)~\right\|_p\leq Ce^{\vert x_0\vert R^\alpha}\|~\f~\|_p,$$
		is fulfilled for each $x_0\in \BR$.
		\item[\textbf (D)] {For every $\g\in L^q(\BR^n;\cl_{0,n})$, with $\displaystyle \frac{1}{p}+\frac{1}{q}=1$, there exists $C>0$ such that
			$$
			\left\vert~\langle \u(\cdot,x_0),\g \rangle_0~\right\vert\leq C2^ne^{\vert x_0\vert R^\alpha}\left\|~\f~\right\|_p \left\|~\g~\right\|_q,
			$$ }
		is fulfilled for each $x_0\in \BR$.
	\end{description}
\end{theorem}

%%% @nelson 
%{First, we establish equivalence between compact support of Fourier transform of $\f$ and Bernstein's type inequality for $\f$.}

The proof of the Theorem \ref{PaleyWienerEquivalences} is divided into two main parts. The first part, concerning the equivalence {\textbf (A)}$\Longleftrightarrow${\textbf (B)}, establishes the connection between the compact support of the Fourier transform of $\f$ and a Bernstein-type inequality for $\displaystyle \f_\pm=\frac{1}{2}(\f\pm \mathcal{H}\f)$. The second part, involving the chain of implications {\textbf (B)}$\Longrightarrow${\textbf (C)}$\Longrightarrow${\textbf (D)}$\Longrightarrow${\textbf (B)}, lays the foundation for the subsequent introduction of a coherent definition of Bernstein-type spaces in the hypercomplex setting.

\textbf{Proof of (A) $\Longrightarrow$  (B)}
First, we observe that, based on Lemma \ref{Uryshon-Lemma}, we can deduce the existence of $\PSI$ such that $\widehat{\PSI}(\XI) = 1$ for all $\XI \in \operatorname{supp}\widehat{\f}$. Next, using the convolution identity, $\f_\pm = {\bm \psi} * \f_\pm$ and referring to the Lemma \ref{DiracHardyRn-Lemma}, we get the following
	\begin{eqnarray*}
		\left(~\D_\theta^\alpha~\right)^k \f_+~=\left(~(-\Delta)^{\frac{\alpha}{2}}~\right)^k{\bm \psi} * \f_+ &\& & \left(~\D_\theta^\alpha~\right)^k \f_-~=e^{-i\pi\theta}\left(~(-\Delta)^{\frac{\alpha}{2}}~\right)^k{\bm \psi} * \f_-,
	\end{eqnarray*}
	for $\displaystyle \f_\pm=\frac{1}{2}(\f\pm \mathcal{H}\f)$,
	so that
	\begin{eqnarray}
		\label{DiracpIneq}\left\|~	\left(~\D_\theta^\alpha~\right)^k \f_\pm~\right\|_p \leq \left\|~	\left(~(-\Delta)^{\frac{\alpha}{2}}~\right)^k \PSI~\right\|_1~\left\|~\f_\pm~\right\|_p
	\end{eqnarray}
	yields by Young's inequality.

	Next, choose $m\in \BN_0$ so that $\left\|~\left(1+\vert \cdot\vert^2\right)^{-m}~\right\|_1<\infty$. Then, from the statement {\textbf (b)} of Theorem \ref{RealPaleyWiener-BumpFunction}, the following set of inequalities is immediate:
	\begin{eqnarray*}
		\left\|~	\left(~(-\Delta)^{\frac{\alpha}{2}}~\right)^k \PSI~\right\|_1 &\leq & \left\|~\left(1+\vert \cdot\vert^2\right)^{-m}~\right\|_1~\sup_{\X \in \BR^n}\left\vert~ \left(1+\vert \X \vert^2\right)^m\left(~(-\Delta)^{\frac{\alpha}{2}}~\right)^k \PSI(\X)~\right\vert_0\\
		&\leq &\Lambda_{m,n}~\left\|~\left(1+\vert \cdot\vert^2\right)^{-m}~\right\|_1~R^{\alpha k}.
	\end{eqnarray*}
	
	So the proof of the statement {\textbf (B)} holds for the constant $$\displaystyle C=\Lambda_{m,n}~\left\|~\left(1+\vert \cdot\vert^2\right)^{-m}~\right\|_1. \qquad \square$$

\textbf{Proof of (B) $\Longrightarrow$ (A)}
	For each $\varepsilon>0$, we assume that $\vert \XI \vert>R+\varepsilon$.
Assuming that the set of estimates $$ \left\|~ \left(~\D_\theta^\alpha~\right)^k \f_\pm~\right\|_p\leq CR^{\alpha k}\left\|~\f_\pm~\right\|_p$$ is satisfied, for each $k\in \BN_0$, we again invoke the Lemma \ref{Uryshon-Lemma} to choose, for $q>1$ such that $\displaystyle \frac{1}{p}+\frac{1}{q}=1$, a real-valued bump function $\PSI\in L^q(\BR^n; \cl_{0,n})$ that satisfies the following conditions:
	%\begin{multicols}{2}
	\begin{description}
		\item[\textbf (i)] $\displaystyle \vert \XI \vert< \sup_{\Y \in \operatorname{supp}\widehat{\f}}\vert \Y \vert$, for all $\XI\in \operatorname{supp}\widehat{\PSI}$;
		\item[\textbf (ii)] $\widehat{\PSI}(\XI)=1$ in $\operatorname{supp}\widehat{\f}$.
	\end{description}
%	\end{multicols}
	
Then, for $\displaystyle \f_\pm = \frac{1}{2}(\f \pm \mathcal{H}\f)$, the convolution identity $\f_\pm = {\bm \psi} * \f_\pm$ allows us to establish, for each $k \in \BN_0$, the set of identities
	\begin{eqnarray*}
		\begin{array}{lll}
			\mathcal{F}\left(\left(~\D_\theta^\alpha~\right)^k\f_+ *{\bm \psi}\right)(\XI)&=&\vert\XI\vert^{\alpha k}\chi_-(\XI)\widehat{\f}(\XI),\\[1ex]
			\mathcal{F}\left(\left(~\D_\theta^\alpha~\right)^k\f_- *{\bm \psi}\right)(\XI)&=&e^{-i\pi \theta}\vert\XI\vert^{\alpha k}\chi_+(\XI)\widehat{\f}(\XI),
		\end{array}
	\end{eqnarray*}
	where $\chi_-(\XI)$ and $\chi_+(\XI)$ denote the idempotents defined by~\eqref{ProjectionOp} (see also~\eqref{ProjectionOpHilbert}).
Then, by applying H\"older's inequality, we derive the following set of inequalities:
	\begin{align*}
		(R+\varepsilon)^{\alpha k} \left\vert~\chi_\mp(\XI) \widehat{\f}(\XI)~\right\vert_0 & \leq \vert \XI\vert^{\alpha k} \left\vert~\chi_\mp(\XI) \widehat{\f}(\XI)~\right\vert_0 
		\leq \left\|~	\left(~\D_\theta^\alpha~\right)^k \f_\pm*{\bm \psi}~\right\|_1 \\
		&\leq  \left\|~	\left(~\D_\theta^\alpha~\right)^k \f_\pm~\right\|_p \left\|~ \PSI~\right\|_q 
		\leq C_\pm R^{\alpha k}\left\|~\f_\pm~\right\|_p \left\|~ \PSI~\right\|_q.
	\end{align*}
	
Therefore, the proof of the statement \textbf{(A)} is a direct consequence of the limit conditions
	\begin{eqnarray*}
		\vert \chi_\mp(\XI)\widehat{\f}(\XI)\vert_0\leq  \lim_{k\rightarrow \infty}C_\pm \left(\frac{R}{R+\varepsilon}\right)^{\alpha k}\left\|~\f_\pm~\right\|_p \left\|~ \PSI~\right\|_q=0. \qquad \square
	\end{eqnarray*}

%%% @nelson
%{\color{purple}Secondly, we prove that Bernstein's inequality is equivalent to the definition of the Bernstein spaces.} 
\textbf{Proof of  (B) $\Longrightarrow$  (C)}
	Let us denote by $P_N$ the operator $$
P_N =\sum_{k=0}^{N}\frac{(-x_0)^k}{k!}\left(~\D_\theta^\alpha~\right)^k.
$$

In the light of the Lemma \ref{DiracHardyRn-Lemma}, it can be deduced that the sequence of inequalities described in \textbf{(B)} is equivalent to
	\begin{eqnarray*}
		~\left\|~\left(~\D_\theta^\alpha~\right)^k\f_\pm~\right\|_p\leq C_\pm R^{\alpha k}\|~\f_\pm~\|_p.
	\end{eqnarray*}
	
	Then the sequence of inequalities
	\begin{align*}
		\left\|~P_N\f_\pm~\right\|_p & \leq \sum_{k=0}^{N}\frac{\vert x_0 \vert^k}{k!}~\left\|~\left(~\D_\theta^\alpha~\right)^k\f_\pm~\right\|_p 
		 \leq \sum_{k=0}^{N}\frac{\vert x_0\vert ^k}{k!}C2^nR^{\alpha k}\left\|~\f_\pm~\right\|_p 
		\leq C2^n~e^{\vert x_0\vert R^\alpha}\left\|~\f_\pm~\right\|_p
	\end{align*}
	follows directly from an inductive application of the Minkowski inequality.
%, as stated in statement (ii) of Proposition \ref{LpNormProp}. 
Therefore, by direct application of the dominated convergence theorem (see \cite[Theorem 1.8]{LiebLoss2001}), we conclude that
	$$\left\|~e^{-x_0\D_\theta^\alpha}\f_\pm~\right\|_p= \lim_{N\rightarrow \infty}\left\|~P_N\f_\pm~\right\|_p\leq C_\pm e^{\vert x_0\vert R}\left\|~\f_\pm~\right\|_p.$$
	
Consequently, for the function $\u$ described by Proposition \ref{CauchyProblemHardy-Convolution} (see \eqref{Fpm} and \eqref{Fpm-Hardy}), the mapping property $\mathcal{H}: L^p(\BR^n;\cl_{0,n})\rightarrow L^p(\BR^n;\cl_{0,n})$ (cf.~\cite[Theorem 2.6]{SteinWeiss72}) ensures the existence of a constant $C > 0$ to satisfy the statement \textbf{(A)}.
$\square$

\textbf{Proof of (C) $\Longrightarrow$ (D)}
Starting with  H\"older's inequality, we obtain, for any $x_0 \in \BR$, the inequality
$$
\left\vert~\left\langle~ \u(\cdot,x_0),\g~\right \rangle_0~\right\vert\leq { 2^n}\left\|~\u(\cdot,x_0)~\right\|_p \left\|~\g~\right\|_q.
$$

Furthermore, assuming that the statement \textbf{(C)} is true, the proof of the statement \textbf{(D)} is immediate.
$\square$

\textbf{Proof of (D) $\Longrightarrow$ (B)}
Recall that the statement ${\textbf (D)}$ asserts that for every $\f\in L^p(\BR^n;\cl_{0,n})$ and $\g\in L^q(\BR^n;\cl_{0,n})$, with $\displaystyle \frac{1}{p}+\frac{1}{q}=1$, the real-valued function
$x_0\mapsto \left\langle \u(\cdot,x_0),\g \right\rangle_0$, equipped with the associated {\textbf inner product} $\left\langle \cdot,\cdot \right\rangle_0$, satisfies
$$\left\langle -\partial_{x_0}\u(\cdot,x_0),\g \right\rangle_0=\left\langle~ \D_\theta^\alpha\u(\cdot,x_0),\g ~\right\rangle_0.$$. 

Then, by induction over $k$, we deduce that
$$(-1)^k(\partial_{x_0})^k\left\langle \u(\cdot,x_0),\g \right\rangle_0=\left\langle \left(~\D_\theta^\alpha~\right)^k\u(\cdot,x_0),\g \right\rangle_0,$$ 
holds for every $k\ in \BN_0$.
	
	Furthermore, we observe that $\u(\cdot,x_0)=e^{-x_0\D_\theta^\alpha}{\bm F}(\cdot,x_0)$ follows from Proposition \ref{CauchyProblemHardy-Convolution}, where ${\bm F}$ denotes the auxiliary function \eqref{Fpm-Hardy}.

	The combination of H\"older's inequality for Clifford-valued functions and the classical Bernstein inequality (cf.~\cite[eq.~(8) on p.~116]{Nik75}) leads to 
	\begin{align*}
		\displaystyle \left\vert ~\left\langle \left(~\D_\theta^\alpha~\right)^k\u(\cdot,x_0),\g \right\rangle_0~\right\vert& = \left\vert~\left\langle (-1)^k(\partial_{x_0})^k \u(\cdot,x_0),\g \right\rangle_0 ~\right\vert    
		\leq  B(x_0) {2^n}  \left\|~\u(\cdot,x_0)~\right\|_p~\left\|~\g~\right\|_q \\
		&\leq  B(x_0) {2^n} R^{\alpha k} \left\|~{\bm F}(\cdot,x_0)~\right\|_p~\left\|~\g~\right\|_q,
	\end{align*}
	where $B$ denotes a continuous function $B:\BR \rightarrow (0,+\infty)$.
	
	Then, by letting $x_0\rightarrow 0^+$ and $x_0\rightarrow 0^-$, respectively, the above inequality reduces to the following set of inequalities
	\begin{eqnarray}
		\label{BernsteinTypeIneq}
		\left\vert~\left\langle \left(~\D^{\alpha}_\theta~\right)^k\f_\pm,\g \right\rangle_0 ~\right\vert\leq C_\pm~ 2^nR^{\alpha k} \left\|~\f_\pm~\right\|_p~\left\|~\g~\right\|_q,
	\end{eqnarray}
	where $\f_\pm=\frac{1}{2}(\f\pm \mathcal{H}\f)$ and $C_\pm$ denote the constants 
	$\displaystyle C_\pm:=\lim_{x_0\rightarrow 0^\pm}B(x_0)$.
	
If $\f={\bm 0}$, the proof of the statement {\textbf (A)} is immediate. Otherwise, for Clifford-valued-functions $\g$ defined by 
\begin{eqnarray*}
	\g(\X)=
	\left(~\D^{\alpha}_\theta~\right)^k\f_\pm(\X)~ \dfrac{\left\vert~\left(~\D^{\alpha}_\theta~\right)^k\f_\pm(\X)~\right\vert^{p-2}_0}{\left\|~\left(~\D^{\alpha}_\theta~\right)^k\f_\pm~\right\|_p^{p-1}}&,~\X\in \BR^n,
\end{eqnarray*}
it holds $\left\|~\g~\right\|_q=
1 $ and
\begin{eqnarray*}
	\left(\left(~\D^{\alpha}_\theta~\right)^k\f_\pm(\X)\right)^\dag \g(\X)= \frac{\left\vert~\left(~\D^{\alpha}_\theta~\right)^k\f_\pm(\X)~\right\vert^{p}_0}{\left\|~\left(~\D^{\alpha}_\theta~\right)^k\f_\pm~\right\|_p^{p-1}} &,~\X\in \BR^n.
\end{eqnarray*}

Thus, by straightforward computation, we establish that \eqref{BernsteinTypeIneq} is equivalent to the set of inequalities provided by the statement {\textbf (B)}, as intended.
$\square$

%\subsubsection*{Proof of $(B)\Rightarrow (C)$}

\section{Bernstein-type spaces and spectral analysis}\label{Applications}

\subsection{Bernstein-type spaces $\mathcal{B}_R^p$}

In Section \ref{sec6}, the employed proof technique for proving Theorem \ref{PaleyWienerEquivalences} allows for the rigorous extension of Bernstein spaces to the hypercomplex framework. More precisely, statement {\textbf (D)} of Theorem \ref{PaleyWienerEquivalences} permits us to define, for every $\g\in \mathcal{S}(\BR^n;\cl_{0,n})$, the
	real-valued function, $\phi_{\g}$, via 
\begin{eqnarray*}
	\displaystyle \phi_{\g}(x_0)=\left\langle  \u(\cdot,x_0),\g \right\rangle_0&,~x_0 \in \BR.
\end{eqnarray*}
Here, $\u$ is the solution to the Cauchy problem \eqref{CauchySpaceFractionalDirac} as established in Proposition~\ref{CauchyProblemHardy-Convolution}. 
This together with the fact that $\langle\cdot,\cdot \rangle_0$ defines an inner product guarantees that $\phi_{\g}\in C^\infty(\BR)$ satisfies the set of relations $\phi_{\g}(0)=\langle \f,\g\rangle_0$ and
\begin{eqnarray*}
(-1)^k\phi_{\g}^{(k)}(x_0)=\left\langle \left(\D_\theta^\alpha\right)^k\u(\cdot,x_0),\g\right\rangle_0, & x_0\in \BR,&k\in \BN_0.
\end{eqnarray*}
By invoking the Riesz representation theorem, it is possible to readily establish that the real-valued function $\phi_{\g}$ is an element of the Bernstein space $B_R^p$ (see Subsection~\ref{StateOfArt}) if and only if Theorem \ref{PaleyWienerEquivalences} holds. This correspondence is closely linked to the approach developed by Pesenson (see references \cite{Pesenson01,Pesenson08,Pesenson15}) and Pesenson-Zayed (see reference \cite{PesensonZayed09}) in their work on abstract Paley-Wiener spaces. Consequently, this framework motivates the introduction of the concept of a Bernstein-type space generated by $\D_\theta^\alpha$, as defined below:

%We begin this section with the definition of the \textit{generalized Bernstein space} generated by $\D_\theta^\alpha$, which serves as a function space that seamlessly describes the null solutions of $\partial_{x_0}+\D_\theta^\alpha$ that are of exponential type $R^\alpha$. 
%The generalization of Bernstein-type spaces begins with the solutions of the Cauchy problem~\eqref{CauchySpaceFractionalDirac}, as outlined in Proposition \ref{CauchyProblemHardy-Convolution}.~This problem is closely related to the approach of Pesenson \cite{Pesenson01,Pesenson08,Pesenson15} and Pesenson-Zayed \cite{PesensonZayed09} on abstract Paley-Wiener spaces.
	
	%%% @nelson
	 %That gives in turn a faithful generalization for the Bernstein spaces, appearing in the proof of \cite[Theorem 5.1]{FranklinHoganLarkin17} and elsewhere.

\begin{definition}[Generalized Bernstein spaces $\mathcal{B}_{R}^{p}$]\label{Bernstein-Spaces}
	For $R>0$ and $1<p<\infty$, the \textit{Bernstein-type space} $\mathcal{B}_{R}^{p}(\D_\theta^\alpha)$ is a Banach right module consisting of all Clifford-valued functions $\u$ such that:
	\begin{description}\itemsep0.5ex
		\item[\textbf (1)] $\u$ is the solution of \eqref{CauchySpaceFractionalDirac} provided by Proposition \ref{CauchyProblemHardy-Convolution};
		\item[\textbf (2)] $\displaystyle \u(\cdot,0)=\f$ belongs to $L^p(\BR^n;\cl_{0,n})$;
		\item[\textbf (3)] {$\displaystyle \sup_{x_0\in \BR}e^{-\vert x_0\vert R^\alpha}\left\|~\u(\cdot,x_0)~\right\|_p~ <\infty$.}
	\end{description}
\end{definition}
Conditions {\textbf (1)} and {\textbf (2)} of Definition \ref{Bernstein-Spaces} are closely related to semigroup theory, while condition {\textbf (3)} is established by statement {\textbf (D)} of Theorem \ref{PaleyWienerEquivalences}. Essentially, the solutions of \eqref{CauchySpaceFractionalDirac} are generated by the action of the semigroup~$\left\{e^{-x_0\D_\theta^\alpha}\right\}_{x_0\in\BR}$ on the boundary values $\displaystyle \f_\pm=\frac{1}{2}\left(\f\pm \mathcal{H}\f\right)$, as demonstrated in Proposition \ref{CauchyProblemHardy-Convolution}. This semigroup plays an analogous role to that of$\left\{e^{iy\frac{d}{dx}}\right\}_{y\in \BR}$, which is induced by the translation operator $f(x)\mapsto f(x+iy)$ in the complex plane. Consequently, the formulation of Bernstein-type spaces, as introduced earlier, follows directly from the equivalences established in Theorem \ref{PaleyWienerEquivalences}.

\subsection{An hypercomplex analogue of the Paley-Wiener theorem}\label{hypercomplexPW}

Let us now turn our attention to the hypercomplex version of the Paley-Wiener theorem. By combining Proposition \ref{RealPaleyWiener-BumpFunction} with Theorem \ref{PW-KouQian-Generalization}, we obtain a generalization of \cite[Theorem 2.1]{KouQian2002} and \cite[Theorem 5.1]{FranklinHoganLarkin17}. In particular, by noting that the condition $\u\in {\mathcal{B}_R^2(\D)}$ is equivalent to statement {\textbf (4)} of \cite[Theorem 5.1]{FranklinHoganLarkin17},
%By merging Proposition \ref{RealPaleyWiener-BumpFunction} and Theorem \ref{PW-KouQian-Generalization}, we are now able to obtain a generalization of \cite[Theorem 2.1]{KouQian2002} and \cite[Theorem 5.1]{FranklinHoganLarkin17}. In particular, 
we can generalize the following chain of equivalences:
\begin{eqnarray*}%\tag{$\mathcal{PW}$}
%	\label{PW_Kou-Qian}	
	\begin{array}{lll}
	\displaystyle	\sup_{(\X,x_0)\in \BR^n \times \BR}e^{-\vert x_0+\X\vert R}\vert \u(\X,x_0)\vert<\infty &\Longleftrightarrow & \operatorname{supp}\widehat{\f}\subseteq \overline{B(0, R)} \\[1ex]
		&\Longleftrightarrow & \u\in \mathcal{B}_R^2(\D).
	\end{array}
\end{eqnarray*}

In this context and beyond, we will use the notation $\vert x_0+\X\vert:=\sqrt{x_0^2+\vert \X\vert^2}$ to represent the norm of the paravector $x_0+\X$ in $\cl_{0,n}$.

	\begin{theorem}\label{PW-KouQian-Generalization}
		Let $\u$ be the solution of the Cauchy problem \eqref{CauchySpaceFractionalDirac} such that 
\begin{eqnarray*}
	\u(\cdot,0)=\f&\mbox{and}&\f\in L^p(\BR^n;\cl_{0,n}). 
\end{eqnarray*}

For $1 < p < \infty$ and $R > 0$, the following statements are equivalent:
	%	\begin{multicols}{2}
			\begin{description}
				\item[\textbf (a)] $\operatorname{supp}\widehat{\f}\subseteq\overline{B(0,R)}$. \\[-1.5ex]
					%\item[~]
				\item[\textbf (b)] $\displaystyle \sup_{(\X,x_0)\in \BR^{n+1}} e^{-\vert x_0+\X\vert R^\alpha}\vert \u(\X,x_0) \vert_0<\infty$. \\
				%\item[~]
				\item[\textbf (c)] $\u \in \mathcal{B}_R^p(\D_\theta^\alpha)$.
			\end{description}
	%	\end{multicols}
	\end{theorem}
	
	\textbf{Proof of (a) $\Longleftrightarrow$ (c)}
	From Proposition \ref{CauchyProblemHardy-Convolution} it easily follows that $\u(\cdot,x_0)=e^{-x_0\D^\alpha_\theta}{\bm F}(\cdot,x_0)$ holds for every $x_0\in \BR$, where ${\bm F}$ is the function defined by \eqref{Fpm}.
		Then, the equivalence ${\textbf (A)}\Longleftrightarrow {\textbf (C)}$ provided by Theorem \ref{PaleyWienerEquivalences}, Sec. \ref{subsec62}, allows us to establish the equivalence
		$$
		\operatorname{supp}\widehat{\f}\subseteq \overline{B(0,R)} \Longleftrightarrow	\displaystyle \sup_{x_0\in\BR}e^{-\vert x_0\vert R^\alpha}\left\|~\u(\cdot,x_0)~\right\|_p\leq C ,
		$$
for some $0 < C < \infty$, thereby proving the equivalence between {\textbf (a)} and {\textbf (c)}.
	$\square$
	
\textbf{Proof of (a) $\Longrightarrow$ (b)}
Let us assume that $\operatorname{supp} \widehat{\f} \subseteq\overline{B(0, R)}$. Then, by Lemma \ref{Uryshon-Lemma}, there exists a real-valued bump function $\PSI$ such that $\widehat{\PSI}(\XI) = 1$ in $\operatorname{supp} \widehat{\f}$. From Proposition \ref{CauchyProblemHardy-Convolution}, we get the identity
$$ 
\u(\cdot,x_0)=
\begin{cases}
\left(e^{-x_0(-\Delta)^{\frac{\alpha}{2}}}\PSI\right)*\f_+&,~x_0>0 \\[1ex]
\left(e^{-x_0e^{-i\pi\theta}(-\Delta)^{\frac{\alpha}{2}}}\PSI\right)*\f_-&,~x_0<0 \\[1ex]
\f&,~x_0=0
\end{cases}
$$
is valid at the level of the distributions, where $\f_\pm=\frac{1}{2}(\f\pm \mathcal{H}\f)$.
As a result, for each $x_0\in \BR \setminus \{0\}$, we get
\begin{eqnarray}
\label{ExpHolder}\left\vert ~\u(\cdot,x_0)~\right\vert_0\leq 
\begin{cases}
	\left\|~e^{-x_0(-\Delta)^{\frac{\alpha}{2}}}\PSI~\right\|_q\left\| \f_+\right\|_p&,~x_0>0 \\ \ \\ %[1ex]
\left\|~e^{-x_0e^{-i\pi\theta}(-\Delta)^{\frac{\alpha}{2}}}\PSI~\right\|_q\left\| \f_-\right\|_p&,~x_0<0,
\end{cases}
\end{eqnarray}
with $1 < p, q < \infty$ such that $\displaystyle \frac{1}{p} + \frac{1}{q} = 1$, it follows by H\"older's inequality.

Next, let us take a close look at the family of norms $\left\|~e^{-x_0e^{-i\pi\beta}(-\Delta)^{\frac{\alpha}{2}}}\PSI~\right\|_q$, with $\beta\in \{0,\theta\}$. By choosing a constant $m\in \BN_0$ such that $\left\|~\left(~1+\vert \X \vert^2~\right)^{-m}~\right\|_q<\infty$ and observing that
\begin{eqnarray*}
	\left\|~e^{-x_0e^{-i\pi\beta}(-\Delta)^{\frac{\alpha}{2}}}\PSI~\right\|_q^q=\\
	=\int_{\BR^n} (1+\vert \X\vert^2)^{-mq}~\left\vert~ (1+\vert\X\vert^2)^me^{-x_0e^{-i\pi\beta}(-\Delta)^{\frac{\alpha}{2}}}\PSI(\X)~\right\vert_0^q \dx,
\end{eqnarray*}
we obtain
\begin{equation}
	\label{BumpExp}\left\|~e^{-x_0e^{-i\pi\theta}(-\Delta)^{\frac{\alpha}{2}}}\PSI~\right\|_q\leq S_\beta \left\|~\left(~1+\vert \X \vert^2~\right)^{-m}~\right\|_q,
\end{equation}
where 
$$
S_\beta=\sup_{\X \in \BR^n} \left\vert~ (1+\vert\X\vert^2)^me^{-x_0e^{-i\pi\beta}(-\Delta)^{\frac{\alpha}{2}}}\PSI(\X)~\right\vert_0.
$$

Furthermore, recalling that the formal Taylor series expansion of $\displaystyle e^{-x_0e^{-i\pi\beta}(-\Delta)^{\frac{\alpha}{2}}}$ is absolutely convergent, it follows that
\begin{eqnarray*}
	S_\beta\leq \sum_{k=0}^\infty \frac{\vert x_0\vert^k}{k!}{\bm s}_k &,~\mbox{with}&{\bm s}_k=\sup_{\X \in \BR^n} \left\vert~ (1+\vert\X\vert^2)^m\left(~(-\Delta)^{\frac{\alpha}{2}}\right)^k\PSI(\X)~\right\vert_0.
\end{eqnarray*}

Following Theorem \ref{RealPaleyWiener-BumpFunction}, it can be deduced that the coefficients ${\mathbf s}_k$ are bounded above by $\Lambda_{m,n}R^{\alpha k}$. As a result, $$ S_\beta\leq \Lambda_{m,n}~e^{\vert x_0\vert R^\alpha}\leq \Lambda_{m,n}~e^{\vert x_0+\X\vert R^\alpha}, $$ for all $(\X,x_0)\in \BR^{n+1}$ and $\beta \in \{0,\theta\}$.

Moreover, the direct combination of the estimates \eqref{BumpExp} and \eqref{ExpHolder}, followed by the mapping property $\mathcal{H}:L^p(\BR^n;\cl_{0,n})\rightarrow L^p(\BR^n;\cl_{0,n})$ (cf.~\cite[Theorem 2.6]{SteinWeiss72}), provides the proof of statement {\textbf (b)}.
$\square$

\textbf{Proof of (b) $\Longrightarrow$ (a)}
Let us choose a real-valued bump function $\PSI \in \mathbb{S}(\BR^n; \cl_{0,n})$ such that $\operatorname{supp} \widehat{\PSI} \subseteq \overline{B(0, R)}$ and $\widehat{\PSI}(\XI) = 1$ in $\overline{B(0, R)}$. Recalling the Lemma \ref{DiracHardyRn-Lemma}, we obtain the following set of identities for each $k \in \BN_0$:
	\begin{eqnarray*}
\left\vert~\left(~\D_\theta^\alpha~\right)^k\u_\pm(\X,x_0)~\right\vert_0=\left\vert~\u_\pm(\X,x_0)*\left(~(-\Delta)^{\frac{\alpha}{2}}~\right)^k\PSI(\X)~\right\vert_0&,~(\X,x_0)\in \BR^{n+1}_\pm.
	\end{eqnarray*}
	
Based on the above, we estimate the following:
	\begin{eqnarray}
		\label{upmEstimate}
		\vert \XI\vert^{\alpha k}\left\vert~\widehat{\u}_\pm(\XI,x_0)~\right\vert_0\leq \left\|~\u_\pm(\cdot,x_0)*\left(~(-\Delta)^{\frac{\alpha}{2}}~\right)^k\PSI~\right\|_1&,~\XI\in \BR^n.
	\end{eqnarray}
	
On the other hand, according to the statement {\textbf (b)}, for every $\sigma>0$ there are two constants $C_+(\sigma)>0$ and $C_-(\sigma)>0$ such that the following set of estimates holds for $\u_\pm$:
\begin{eqnarray}
	\label{GrowthCondition-sigma}
\left\vert~\u_\pm(\X,x_0)~\right\vert_0 \leq C_\pm(\sigma)e^{\sigma R^\alpha}&,~0<\vert x_0+\X \vert\leq \sigma.
\end{eqnarray}

Also, since the real-valued bump function $\PSI$ satisfies the conditions of Proposition \ref{RealPaleyWiener-BumpFunction}, choosing $m \in \BN_0$ such that $\left\|~\left(~1+\vert~\cdot~\vert^2~\right)^{-m}~\right\|_1<\infty$, yields the following $L^1-L^1$ estimate
\begin{eqnarray}
\label{L1-estimate} \left\|~\left(~(-\Delta)^{\frac{\alpha}{2}}~\right)^k\PSI~\right\|_1\leq \Lambda_{m,n}R^{\alpha k}~\left\|~\left(~1+\vert~\cdot~\vert^2~\right)^{-m}~\right\|_1.
\end{eqnarray}

Thus, the direct application of Young's inequality to the right hand side of \eqref{upmEstimate}, followed by the estimates of \eqref{GrowthCondition-sigma} and \eqref{L1-estimate}, yields 
		\begin{eqnarray}
		\label{upmEstimate-Young}
		\vert \XI\vert^{\alpha k}\left\vert~\widehat{\u}_\pm(\XI,x_0)~\right\vert_0\leq \vartheta_{m,n}^\pm(\sigma) R^{\alpha k}e^{\sigma R^\alpha}&,~\XI\in \BR^n &,~0<\vert x_0 \vert \leq \sigma.
	\end{eqnarray}
where $\vartheta_{m,n}^+$ and $\vartheta_{m,n}^-$ are continuous functions $\vartheta_{m,n}^\pm:\BR\rightarrow (0,+\infty)$ defined by
		$$\vartheta_{m,n}^\pm(\sigma)=C_\pm(\sigma)\Lambda_{m,n}\left\|~\left(~1+\vert~\cdot~\vert^2~\right)^{-m}~\right\|_1.
	$$
	
To establish the statement {\textbf (a)} of Theorem \ref{PW-KouQian-Generalization}, let us consider the condition $\vert \XI\vert > R + \varepsilon$ for every $\varepsilon > 0$.
Then, for sufficiently large values of $k \in \BN_0$ and $\sigma>0$, we obtain from \eqref{upmEstimate-Young} the estimate
			\begin{eqnarray}
			\label{upmEstimate-epsilon}
			\left\vert~\widehat{\u}_\pm(\XI,x_0)~\right\vert_0\leq \vartheta_{m,n}^\pm(x_0) \left(\frac{R}{R+\varepsilon}\right)^{\alpha k}e^{\sigma R^\alpha}&,~\vert \XI\vert>R+\varepsilon &,~0<\vert x_0 \vert \leq \sigma.
		\end{eqnarray}
		
So, by letting $k \rightarrow \infty$ on both sides of \eqref{upmEstimate-epsilon}, we get
\begin{eqnarray*}
	\widehat{\u}_\pm(\XI,x_0)=0,& \mbox{for}&\vert \XI\vert>R+\varepsilon,~~\mbox{and}~~\pm x_0>0.
\end{eqnarray*}
		
Based on Proposition \ref{CauchyProblemHardy}, we deduce that $\f\in L^p(\BR^n;\cl_{0,n})$ is given by
$$
\f=\lim_{x_0\rightarrow 0^+}{\u}_+(\cdot,x_0)+\lim_{x_0\rightarrow 0^-}{\u}_-(\cdot,x_0).
$$

So we find that $\widehat{\f}(\XI)=0$ for $\vert \XI\vert>R+\varepsilon$. This concludes our proof.
$\square$

\begin{remark}
Despite the initial obstacles highlighted in \cite{KouQian2002} in proving the support condition $\operatorname{supp}\widehat{\f}\subseteq \overline{B(0,R)}$ for Theorem \ref{PW-KouQian-Generalization} due to the absence of a Cauchy integral-type formula in the space-fractional setting, a significant extension of \cite[Theorem 2.1]{KouQian2002} was successfully achieved. This extension was carried out by adopting an approach similar to the proof of ${\textbf (B)}\Longrightarrow {\textbf (A)}$ in Theorem \ref{PaleyWienerEquivalences}, rather than relying on the Taylor series expansion of the generalized Cauchy kernels $\E_{\alpha,n}^\pm(\Y,z)$ introduced in Subsection \ref{subsec43}.
\end{remark}

\subsection{Spectral analysis}\label{SpectralAnalysis}

Next, we explicitly determine the maximal radius $R$ for which $\operatorname{supp} \widehat{\f}\subseteq \overline{B(0, R)}$. This determination is achieved by generalizing the Landau-Kolmogorov-Stein inequality that underlies the Favard constants (cf.~\cite{Stein57}). The latter are defined as follows:
$$
K_j=\frac{4}{\pi}\sum_{r=0}^{\infty} \frac{(-1)^{r(j+1)}}{(2r+1)^{j+1}}~(j\in \BN).
$$

Therefore, an explicit formula for the maximal radius can be easily derived, analogous to the real Paley-Wiener theory (cf.~\cite[Theorem 2.6]{AndersenJeu10}).  This result follows from the analysis of the growth of the sequence of Clifford-valued functions $\left(~(\D_\theta^\alpha)^k\f\right)_{k\in \BN_0}$.

 The following result, a hypercomplex analog of the result of \cite[Lemma 6]{PesensonZayed09}, serves as the basis for the following analysis.

\begin{lemma}[Landau-Kolmogorov-Stein type inequality]\label{SteinKolmogorov-Lemma}
	For each $\f\in L^p(\BR^n;\cl_{0,n})$, one has the following inequality
	\begin{eqnarray}
\label{SteinKolmogorov-ineq}		\left\|~\left(~\D^{\alpha}_\theta~\right)^k\f~\right\|_p^{\ell} \leq C_{k,\ell}~\|~\f~\|_p^{\ell-k} ~\left\|~\left(~\D_\theta^\alpha~\right)^\ell\f~\right\|_p^{k}, & 0\leq k\leq \ell,
	\end{eqnarray}
where $C_{k,\ell}$ represents the constant $\dfrac{\left(K_{\ell-k}\right)^\ell}{\left(K_{\ell}\right)^{\ell-k}}$.
\end{lemma}

Proof:
	Starting from Landau-Kolmogorov-Stein's inequality (cf.~\cite[PART II]{Stein57}), we observe that
	$$\left\|~(\partial_{x_0})^k\phi~\right\|_\infty^{\ell} \leq C_{k,\ell}~\|~\phi~\|_\infty^{\ell-k} ~\left\|~(\partial_{x_0})^\ell \phi~\right\|_\infty^{k}$$
	yields for the analytic type function $\phi$, defined via 
	\begin{eqnarray*}
	\displaystyle \phi(x_0)=\left\langle  \u(\cdot,x_0),\g \right\rangle_0&,~x_0 \in \BR,
	\end{eqnarray*}
	where $\u$ is the solution of the Cauchy problem \eqref{CauchySpaceFractionalDirac} given by the Proposition \ref{CauchyProblemHardy-Convolution}.
	This is equivalent to 
	\begin{eqnarray*}\left\|~\left\langle \left(~\D^{\alpha}_\theta~\right)^k\left(~\u(\cdot,x_0)~\right),\g \right\rangle_0~\right\|_\infty^{\ell} &\leq& C_{k,\ell}~\left\|~\left\langle \u(\cdot,x_0),\g \right\rangle_0~\right\|_\infty^{\ell-k}  \\
		&\times& ~\left\|~\left\langle \left(~\D^{\alpha}_\theta~\right)^\ell\left(~\u(\cdot,x_0)~\right),\g \right\rangle_0~\right\|_\infty^{k} .
	\end{eqnarray*}
	
In particular, for $x_0=0$ the direct application of the H\"older inequality leads to
	\begin{eqnarray*}\left\vert~\left\langle \left(~\D^{\alpha}_\theta~\right)^k\f,\g \right\rangle_0~\right\vert^{\ell} &\leq& C_{k,\ell}~\left\|~\f~ \right\|_p^{\ell-k} \left\|~\left(~\D^{\alpha}_\theta~\right)^\ell\f~\right\|_p^{k} \left\|~\g~\right\|_{q}^{\ell}.
	\end{eqnarray*}

If $\f={\bm 0}$, the inequality \eqref{SteinKolmogorov-ineq} is automatically satisfied. Otherwise,  $$\g(\X)=
\left(~\D^{\alpha}_\theta~\right)^k\f(\X)~ \dfrac{\left\vert~\left(~\D^{\alpha}_\theta~\right)^k\f(\X)~\right\vert^{p-2}_0}{\left\|~\left(~\D^{\alpha}_\theta~\right)^\ell\f~\right\|_p^{p-1}}$$ we also get \eqref{SteinKolmogorov-ineq}. This proves our claim.
$\square$

After establishing the essential components, the remaining theorem of this subsection is proved, following the proof structure in \cite[Theorem 3.4]{Pesenson08}, \cite[Theorem 8]{PesensonZayed09}, and \cite[Theorem 2.9]{Pesenson15} in the Hilbert space context.
\begin{theorem}\label{SpectralRadiusTheorem}
	For $1<p<\infty$, let $\f\in L^p(\BR^n;\cl_{0,n})$ such that $\operatorname{supp}\widehat{\f}\subseteq \overline{B(0,R)}$, for some $0<R<\infty$.
	Then, the limit
	\begin{equation}%\tag{$\mathcal{L}_{\alpha,k}$}
		\label{Limit-alphak}\lim_{k\rightarrow \infty} {\left\|~\left(~\D^{\alpha}_\theta~\right)^k\f~\right\|_p^{\frac{1}{\alpha k}}},
	\end{equation}
	exists, is finite and coincides with
	\begin{equation}%\tag{$\mathcal{R}$}
		\label{Rf}R(\widehat{\f}):=\inf\left\{~\sigma>0~:~\operatorname{supp}\widehat{\f}\subseteq \overline{B(0,\sigma)}~\right\}.
	\end{equation}
	
	Conversely, if $\u\in \mathcal{B}^p(\D^\alpha_\theta)$ is such that $\u(\cdot, 0) = \f$, and the limit \eqref{Limit-alphak} exists and is finite, then $$\operatorname{supp}\widehat{\f}\subseteq \overline{B\left(~0,R(\widehat{\f})~\right)}.$$
\end{theorem}

Proof:
	Starting from Lemma \ref{SteinKolmogorov-Lemma}, the inequality \eqref{SteinKolmogorov-ineq}	
	leads to
	$$\left\|~\left(~\D^{\alpha}_\theta~\right)^k\f~\right\|_p^{\frac{1}{\alpha k}}=\left(\left\|~\left(~\D^{\alpha}_\theta~\right)^k\f~\right\|_p^\ell\right)^{\frac{1}{\alpha k\ell}} \leq (C_{k,\ell})^{\frac{1}{\alpha k\ell}}~\|~\f~\|_p^{\frac{1}{\alpha k}-\frac{1}{\alpha \ell}} ~\left\|~\left(~\D^{\alpha}_\theta~\right)^\ell\f~\right\|_p^{\frac{1}{\alpha \ell}}.$$
	
	This implies that
	$$\left\|~\left(~\D^{\alpha}_\theta~\right)^k\f~\right\|_p^{\frac{1}{\alpha k}}\leq ~\|~\f~\|_p^{\frac{1}{\alpha k}} ~\liminf_{\ell \rightarrow \infty}\left\|~\left(~\D^{\alpha}_\theta~\right)^\ell\f~\right\|_p^{\frac{1}{\alpha \ell}}.$$
	
Furthermore, if we take $\displaystyle \limsup_{k\rightarrow \infty}$ on both sides of the above inequality, we get
	$$\limsup_{k\rightarrow \infty}\left\|~\left(~\D^{\alpha}_\theta~\right)^k\f~\right\|_p^{\frac{1}{\alpha k}}\leq~\liminf_{\ell \rightarrow \infty}\left\|~\left(~\D^{\alpha}_\theta~\right)^\ell\f~\right\|_p^{\frac{1}{\alpha \ell}},$$
	which shows that the limit \eqref{Limit-alphak} exists.
	
Next, we prove that $R(\widehat{\f})$, defined by \eqref{Rf}, is an upper bound of the limit \eqref{Limit-alphak}. To do this, we observe that, given the statement {\textbf (B)} of Theorem \ref{PaleyWienerEquivalences}, one has the inequality
	$$
	\left\|~\left(~\D^{\alpha}_\theta~\right)^k\f~\right\|_p^{\frac{1}{\alpha k}}\leq C^{\frac{1}{k}}R(\widehat{\f})~\|~\f~\|_p^{\frac{1}{\alpha k}},
	$$
	for some $0<C<\infty$.
	
	Then, by taking the limit $k\rightarrow \infty$ of both sides of the above inequality, we conclude that
	$$\displaystyle \lim_{k\rightarrow \infty} {\left\|~\left(~\D^{\alpha}_\theta~\right)^k\f~\right\|_p^{\frac{1}{\alpha k}}}\leq R(\widehat{\f}).$$
	
Now, to prove that the limit \eqref{Limit-alphak} coincides with \eqref{Rf}, we assume by contradiction that
 $$\displaystyle \lim_{k\rightarrow \infty} {\left\|~\left(~\D^{\alpha}_\theta~\right)^k\f~\right\|_p^{\frac{1}{\alpha k}}}< R(\widehat{\f}).$$
	
 Then there exist $\sigma>0$ and $M>0$ such that 
 \begin{eqnarray*}
 \displaystyle \displaystyle \lim_{k\rightarrow \infty} {\left\|~\left(~\D^{\alpha}_\theta~\right)^k\f~\right\|_p^{\frac{1}{\alpha k}}}<\sigma<R(\widehat{\f}) &\&& \left\|~\left(~\D^{\alpha}_\theta~\right)^k\f~\right\|_p\leq M\sigma^{\alpha k} \| \f\|_p.
 \end{eqnarray*}
	
	So by the statement {\textbf (A)} of the Theorem \ref{PaleyWienerEquivalences} we get $\operatorname{supp}\widehat{\f}\subseteq \overline{B(0,\sigma)}$, which contradicts the definition \eqref{Rf} of $R(\widehat{\f})$.
	
Conversely, let us assume that the limit \eqref{Limit-alphak} exists and is finite. This means that
	$$\limsup_{k\rightarrow \infty}\left\|~\left(~\D^{\alpha}_\theta~\right)^k\f~\right\|_p^{\frac{1}{\alpha k}}=~\liminf_{k \rightarrow \infty}\left\|~\left(~\D^{\alpha}_\theta~\right)^k\f~\right\|_p^{\frac{1}{\alpha k}}.$$
	
	Then there exists $R>0$ such that the Bernstein type inequality, which appears in statement {\textbf (B)} of Theorem \ref{PaleyWienerEquivalences}, holds. Thus, according to statement {\textbf (A)} of Theorem \ref{PaleyWienerEquivalences}, this is equivalent to saying that $\operatorname{supp}\widehat{\f}\subseteq \overline{B(0, R)}$.
	
	Then, from the definition \eqref{Rf} of $R(\widehat{\f})$, it immediately~follows that the inequality $$
	~\left\|~\left(~\D^{\alpha}_\theta~\right)^k\f~\right\|_p\leq C\left(~R(\widehat{\f})~\right)^{\alpha k}\|~\f~\|_p ,$$
	holds for every $k \in \BN_0$, is equivalent to
	$$\dfrac{\left\|~\left(~\D^{\alpha}_\theta~\right)^k\f~\right\|_p^{\frac{1}{\alpha k}}}{C^{\frac{1}{\alpha k}}\left\|~\f~\right\|_p^{\frac{1}{\alpha k}}} \leq R(\widehat{\f}).
	$$
	
Furthermore, by allowing $k\rightarrow\infty$ on both sides of the above inequality, we get
	$$
	\lim_{k\rightarrow \infty}{\left\|~\left(~\D^{\alpha}_\theta~\right)^k\f~\right\|_p^{\frac{1}{\alpha k}}}\leq R(\widehat{\f}).
	$$
	
	Finally, the proof that the left side of the above inequality is indeed equal to $R(\widehat{\f})$ is obtained by applying the contradiction argument outlined previously.  
$\square$

\section{Some Outlook}\label{Conclusion}

The foundations of a Paley-Wiener theory have been developed, and this theory accurately describes Clifford-valued functions of exponential type $R^\alpha$, $\u$, from the knowledge of the support of the Fourier transform of $\displaystyle \f=\u(\cdot,0)$, and vice versa. In particular, the spectral formula for the maximum radius $R$, provided by the Theorem \ref{SpectralRadiusTheorem}, allows for an explicit calculation of the so-called bandwidth of the generalized Bernstein space $\mathcal{B}_R^p(\D_{\theta}^\alpha)$.

The predominant comprehension of  $\mathcal{B}_R^p(\D_{\theta}^\alpha)$ indicates the possibility for a hypercomplex extension of the Whittaker-Shannon sampling theorem, analogous to the one outlined in \cite{KouQianSommen07}. However, further investigation is necessary to develop numerical implementations of sampling reconstruction schemes. This investigation necessitates a comprehensive understanding of {\textit Reproducing Kernel Hilbert Spaces}, which is beyond the scope of the present paper. The authors intend to explore these topics in a future publication.

In a different direction, it would be worthwhile to consider the potential applications of our results for values of $0 < p \leq 1$ and $p = \infty$. While many steps in our proofs appear to necessitate only minor alterations, as evidenced by the results obtained in \cite{DangMaiQian2020}, this approach is contingent on intricate properties involving duality theorems between Hardy-type spaces $H^1$ and $BMO$-type spaces, as well as, more generally, between Hardy-type spaces $H^p$ and Morrey-Campanato type spaces (see, for example, \cite{DuongYan05,Yan08} for a comprehensive overview).

\section*{Acknowledgements}

N. Faustino was supported by CIDMA under the FCT Multi-Annual Financing Program for R\&D Units.

\end{document}